\documentclass{article}
\usepackage{amsmath,amsfonts,epsfig}
\newcommand{\bb}{\mathbb}

\newcommand{\cx}{{\bb C}}
\newcommand{\half}{{\bb H}}
\newcommand{\integers}{{\bb Z}}
\newcommand{\natls}{{\bb N}}

\newcommand{\reals}{{\bb R}}

\newcommand{\hthree}{{\bb H}^3}
\newcommand{\htwo}{{\bb H}^2}

\newcommand{\tube}{{\mathbb T}}

\renewcommand{\bold}[1]{\medskip \noindent {\bf \boldmath #1
                        }\nopagebreak[4]}

\newcommand{\qed}[1]{\nopagebreak[4]{\tiny \hfill
\fbox{\ref{#1}} \linebreak }\pagebreak[2]}
\newcommand{\bdry}{\partial}

\newcommand{\compos}{\circ}

\newcommand{\del}{\partial}

\newcommand{\disjunion}{\sqcup}

\newcommand{\nullset}{\emptyset}

\newcommand{\st}{\; : \;}         %Such that

\newcommand{\chat}{\widehat{\cx}}

\newcommand{\zed}{\integers}
\newcommand{\area}{\operatorname{area}}

\newcommand{\core}{\operatorname{core}}

\newcommand{\diam}{\operatorname{diam}}

\newcommand{\inj}{\operatorname{inj}}

\newcommand{\Isom}{\operatorname{Isom}}

\newcommand{\Mod}{\operatorname{Mod}}

\newcommand{\PSL}{\operatorname{PSL}}

\newcommand{\sym}{\operatorname{sym}}

\newcommand{\Teich}{\operatorname{Teich}}

\newcommand{\vol}{\operatorname{vol}}

\newtheorem{theorem}{Theorem}[section]
\newtheorem{prop}[theorem]{Proposition}
\newtheorem{lemma}[theorem]{Lemma}
\newtheorem{cor}[theorem]{Corollary}

\newenvironment{refthm}[1]{{\bold{#1}}\em}{{\bigskip}}

\newcommand{\cC}{{\cal C}}

\newcommand{\cK}{{\cal K}}
\newcommand{\cL}{{\cal L}}
\newcommand{\cM}{{\cal M}}
\newcommand{\cN}{{\cal N}}

\newcommand{\cP}{{\cal P}}

\newcommand{\pml}{{\cP \cM \cL}}

\newcommand{\calC}{{\mathcal C}}

\newcommand{\calK}{{\mathcal K}}

\newcommand{\calM}{{\mathcal M}}
\newcommand{\calN}{{\mathcal N}}

\newcommand{\calP}{{\mathcal P}}

\newcommand{\calS}{{\mathcal S}}

\renewcommand{\hbar}{\bar{{\mathbb H}}^3}

\renewcommand{\core}{C}
\newcommand{\bilip}{\operatorname{bilip}}

\renewcommand\marginpar[1]{} % Kill marginpars for the final version

\begin{document}

\title{\bf Geometric inflexibility and 3-manifolds that fiber over the
circle}
\author{J. Brock\thanks{Research supported by the John S. Guggenheim
    foundation and by the NSF.} \  and K. Bromberg\thanks{Research supported by the
    Alfred P. Sloan Foundation and by the NSF.}}

\date{August 11, 2010}

\maketitle

\begin{abstract}
We prove hyperbolic 3-manifolds are geometrically inflexible: a unit
quasiconformal deformation of a Kleinian group extends to an
equivariant bi-Lipschitz diffeomorphism between quotients whose
pointwise bi-Lipschitz constant decays exponentially in the distance
form the boundary of the convex core for points in the thick part.
Estimates at points in the thin part are controlled by similar estimates on the
complex lengths of short curves.  We use this inflexibility to give a
new proof of the convergence of pseudo-Anosov double-iteration on the
quasi-Fuchsian space of a closed surface, and the resulting
hyperbolization theorem for closed 3-manifolds that fiber over the circle
with pseudo-Anosov monodromy.

\end{abstract}

\section{Introduction}
\label{intro}

In the study of hyperbolic structures on $3$-manifolds, the rigidity
theorems of Mostow and Sullivan allow for coarse methods to play a key
role in the classification of structures up to isometry: it suffices
to exhibit a uniformly bi-Lipschitz map between two hyperbolic
3-manifolds with the same asymptotic data to conclude they are in fact
isometric.

A general theme in work of Thurston has
been the notion of limiting to rigidity, wherein a family of
hyperbolic structures has a quasiconformally rigid limit.  Such
discussions suggest a qualitative notion of {\em inflexibility } for
manifolds far out in the sequence: a unit quasiconformal deformation
at infinity has exponentially deteriorating effect at the basepoint as the geometry
freezes around it.

This qualitative notion was made more precise for manifolds with
injectivity radius bounds by McMullen (see \cite{McMullen:book:RTM}),
but the assumption of injectivity bounds is very restrictive.  Though
upper bounds on the injectivity radius in the convex core follow from
tameness (now known for arbitrarily $M$ with finitely generated
$\pi_1$ \cite{Agol:tame,Calegari:Gabai:tame}), the lower bound is
non-generic
\cite{McMullen:cusps,Canary:Culler:Hersonsky:Shalen,Canary:Hersonsky:cusps}.
In this paper we prove an exponential decay theorem for the $L^2$-norm
of a harmonic deformation a hyperbolic 3-manifold. This allows us to
prove inflexibility theorems for arbitrary hyperbolic
3-manifolds. Here is a sample theorem which generalizes McMullen's
result.
\begin{theorem}{\sc (Geometric Inflexibility)}
  Given a hyperbolic 3-manifold $M$, a $K$-bi-Lipschitz diffeomorphic
  hyperbolic 3-manifold $M'$, and an $\epsilon>0$, there is a
  diffeomorphism $\Phi \colon M \to M'$ whose bi-Lipschitz distortion
  in the $\epsilon$-thick part of the convex core $C(M)$ decays
  exponentially with the distance from $\bdry C(M)$ with rate of decay
  depending only on $\epsilon$, $K$ and the topology of $\bdry M$.
\label{introinflex}
\end{theorem}
See Theorem~\ref{thickinflex} for a more precise version.  

Although Theorem~\ref{introinflex} does not give estimates on the
bi-Lipschitz constant in the thin part, this is to be expected.
Indeed, there are harmonic deformations whose distortion within a
Margulis tube is roughly constant over the tube and does not decay in
the depth into the tube -- the pointwise bounds on the distortion (the
{\em strain}) are determined by its behavior on the boundary of the
tube.  In this sense, Theorem~\ref{introinflex} is sharp, and in
fact optimal, in that we can only expect at best exponential decay of
the bi-Lipschitz constant in the thick part.

On the other hand, the proof of Theorem~\ref{introinflex} is quite
robust and applies to a variety of other situations. For example we
can control the ratio of the change in length of moderate length
geodesics by constants that exponentially decay in the depth of the
geodesic in the convex core. We obtain similar control over short
geodesics by measuring the depth of their entire Margulis tubes. In a future paper we will apply our
methods to deformations of hyperbolic cone-manifolds where
the depth is measured by distance from the singular locus. For both
smooth, complete hyperbolic manifolds and for cone-manifolds, the
Schwarzian derivative can be similarly controlled at components of the
conformal boundary that are fixed under the deformation.

We emphasize that while McMullen's inflexibility theorem is ultimately
a consequence of the compactness of hyperbolic 3-manifolds with
injectivity radius bounds and basepoints in the convex core, our
arguments harness explicit analytic estimates on the pointwise
$L^2$-norm of the deformation to
obtain sharp estimates on the bi-Lipschitz distortion of a deformation
at infinity.

\bold{Convergence results.}  Inflexibility provides for new approaches
and techniques in the theory of Kleinian groups.  To outline these
results, we briefly recall notions from their deformation theory.

Given a closed surface $S$ of negative Euler-characteristic the
Teichm\"uller space, $\Teich(S)$, parametrizes pairs $(f,X)$ of marked
hyperbolic surfaces $$f \colon S \to X,$$ where $f$ is a homeomorphism
up to marking preserving isometry.  The modular group $\Mod(S)$ of
isotopy classes of orientation preserving self homeomorphisms of $S$
acts naturally on $\Teich(S)$ by $\varphi (f,X) = (f \circ
\varphi^{-1},X)$.  A mapping class is {\em pseudo-Anosov} if for each
essential isotopy class of simple closed curves $\gamma$ we have
$\varphi^n(\gamma) \not\simeq \gamma$ for $n \not= 0$.

L. Bers proved that for each pair $(X,Y) \in \Teich(S) \times \Teich(S)$
there is a unique {\em quasi-Fuchsian simultaneous uniformization}, namely, a single
Kleinian group $\Gamma \cong \pi_1(S)$ for which $\Gamma$ leaves
invariant a directed Jordan curve $\Lambda$ in $\chat$  with the
property that 
$\chat \setminus \Lambda = \Omega_X \disjunion \Omega_Y$, where 
$\Omega_X /\Gamma = X$ and $\Omega_Y /\Gamma = Y$ (see
\cite{Bers:simunif}).

As a tool in the deformation theory of
Kleinian groups, Theorem~\ref{introinflex} guarantees convergence in
certain cases where the depth in the convex core at the basepoint
diverges quickly enough.  In particular, Theorem~\ref{introinflex}
gives a new proof of Thurston's double limit theorem for pseudo-Anosov
iteration, the main step in the hyperbolization for 3-manifolds that
fiber over the circle with pseudo-Anosov monodromy (see
\cite{Thurston:hype2,Otal:book:fibered,McMullen:book:RTM}).

\begin{theorem}{ \sc (Pseudo-Anosov Double Limits)}
  For each $X$ and $Y$ in the Teichm\"uller space $\Teich(S)$, and
  each pseudo-Anosov mapping class $\psi \in \Mod(S)$, the double
  iteration $Q(\psi^{-n} (X), \psi^n (Y) )$ converges algebraically
  and geometrically to a limit $Q_\infty \in AH(S)$.
\label{introdoubleconverges}
\end{theorem}
See Theorem~\ref{doubleconverges}.  Note that the convergence {\em up
  to subsequence} was proven earlier by Thurston (see
\cite{Thurston:hype2}).  Convergence was later proven in
\cite{Cannon:Thurston:peano}; McMullen gave a more explicit treatment
 in \cite{McMullen:book:RTM}.  
Note that in our result the
quasi-conformal rigidity of the limit is a direct consequence of the
geometric inflexibility theorem.

Because for each $n$ the manifold $Q_n = Q(\psi^{-n}(X),\psi^n(Y))$ admits a
uniformly bi-Lipschitz diffeomorphism $\Psi_n$ in the homotopy class of $\psi$,
we may apply the inflexibility theorem 
to obtain an isometry $\Psi \colon Q_\infty \to
Q_\infty$ in the homotopy class of $\psi$.  The quotient
$Q_\infty/\langle \Psi \rangle$ is a hyperbolic 3-manifold with the
homotopy type of $T_\psi$, which is thus homeomorphic to $T_\psi$ by a
theorem of Stallings.  We arrive at Thurston's original theorem.
\begin{theorem}[Thurston]{\sc (Mapping Torus Hyperbolic)}
  Let $\psi \in \Mod(S)$ be pseudo-Anosov.  Then the mapping torus
  $T_\psi = S \times [0,1]/(x,0) \sim (\psi(x),1)$ admits a complete
  hyperbolic structure.
\label{hyperbolization}
\end{theorem}

\bold{Curve complex distance and convex core width.}  To describe how
Theorem~\ref{introdoubleconverges} follows from
Theorem~\ref{introinflex}, we remark that one key step is show
linear growth of the width of the convex core in terms of the iterate of
the pseudo-Anosov applied to each factor.  As the width of the core
grows, the geometric effect of the next iterate decays at the basepoint
exponentially fast, and convergence follows.

To show the growth in width is linear, however, the combinatorial
properties of curves on surfaces play a crucial role.  The collection
of isotopy classes $\calS$ of essential simple closed curves on $S$
can be encoded as a graph $\calC(S)$ with vertices corresponding to
elements of $\calS$ and edges joining vertices if their corresponding
classes can be represented by disjoint curves on $S$.  This graph has
the structure of a $\delta$-hyperbolic metric space if each edge is
assigned length $1$ \cite{Masur:Minsky:CCI}.  Though $\calC(S)$ can be
given the structure of a complex by associating $k$-simplices to
$k+1$-tuples of vertices whose representatives can be realized
disjointly, these higher dimensional simplices do not play a role
here.

Among the many reflections of the combinatorics of 
$\calC(S)$ in the geometry of hyperbolic 3-manifolds, the {\em width} of the
convex core of a quasi-Fuchsian manifold is an important new example.
We show the following.
\begin{cor}{\sc (Wide Cores)}
\label{introquasifuchsianthick}
Given a closed surface $S$, there is linear function $f$ such that the
distance between the boundary components of the convex core
$C(Q(X,Y))$ of a quasi-Fuchsian manifold $Q(X,Y)$ in $QF(S)$ is bounded
below by $f(d_{\cC}(X,Y))$.
\end{cor}
(See Corollary~\ref{quasifuchsianthick}).  Here, the distance
$d_\calC(X,Y)$ is shorthand: if $S$ has genus $g$, there is a uniform
$L_g >0$ so that for each $X \in \Teich(S)$ has a the length of the
shortest essential closed loop on $X$ is bounded by $L_g$.
Furthermore, any two shortest loops have uniformly bounded
intersection, by the collar lemma.  It follows that there is a
coarsely defined map from $\Teich(S)$ to the complex of curves, that
sends each $X$ to the collection of vertices whose simple closed
curves have length less than $L_g$ on $X$.  Then $d_\calC(X,Y)$
measures the maximal distance in $\calC(S)$ between shortest curves on
$X$ and on $Y$.

Since the action of pseudo-Anosov iteration has linear
growth in the curve complex, it follows that the width of the convex
core of the double pseudo-Anosov iteration
$$Q(\psi^{-n} (X), \psi^n (Y))$$ is linear in $n$.  Combining these
estimates on core width with Theorem~\ref{introinflex}, Geometric
Inflexibility,  we obtain Thurston's original result.  

It should be noted, however, that Theorem~\ref{introdoubleconverges}
is a convergence theorem rather than a compactness theorem.  In
particular, the rigidity of the limit is implicit in the proof.  As
such, where Thurston's original proof appealed to Sullivan's rigidity
theorem after showing the limit has limit set all of $\chat$, the
existence of a hyperbolic structure on the mapping torus for $\psi$
here is self-contained.

We remark that the linear growth in the width of the convex core with
distance between the bounded length curves on its boundary in
$\calC(S)$ is not specific to pseudo-Anosov deformations.  In
particular, the methods of Theorem~\ref{introdoubleconverges} extend
immediately to apply to sequences $\{Q(X_n,Y_n)\}_n$ of quasi-Fuchsian
manifolds for which we have the bounds $d_T(X_n,X_{n+1}) \le K$ and
$d_T(Y_n,Y_{n+1}) \le K$, and the curve complex distance
$d_\calC(X_n,Y_n)$ grows linearly with $n$.

\bold{Ending laminations and efficient approximations.}  We remark
that a key further application of Theorem~\ref{introinflex} will be a
new approach to the {\em ending lamination conjecture}
\cite{Brock:Canary:Minsky:elc} via efficient approximations by maximal
cusps.  In short, Minsky's {\em a priori bounds} theorem
\cite{Minsky:CKGI} guarantees that for any hyperbolic 3-manifold $M$
in the boundary of a {\em Bers slice} $B_Y = \{Q(X,Y) \st Y \in
\Teich(S)\}$ 
%with the ending lamination $\lambda \in \bdry \calC(S)$
there is an essentially canonical sequence of maximal simplices $P_n
\in \calC(S)$ with $P_n \to \lambda$, $\lambda \in \bdry \calC(S)$,
(the boundary point $\lambda$ is the {\em ending lamination} for $M$)
whose corresponding curves arise with uniformly bounded length
$\ell_M(P_n) < L$ in $M$.

By an application of the grafting technique of
\cite{Bromberg:bers,Brock:Bromberg:density} together with a covering
argument as in \cite{Bromberg:Souto:density} we may, in effect, drill
$P_n$ out of $M$ to obtain a maximal cusp $C_n \in \bdry B_Y$, by a
deformation that has a bounded effect on the geometry in a compact
core $\calM \subset M$.  By the inflexibility theorem, the effect of
this process on the geometry of $\calM$ decays with the distance of
the geodesic representatives of the curves in $P_n$ from $\calM$.  It
follows that he sequence $C_n$ converges back to $M$.  Since $P_n$
depend only on $\lambda$, it follows that $\lambda$ determines $M$.
We take up this approach in \cite{BBES:elc}.

\bold{Plan of the paper.}  A significant component of the paper
involves the study of harmonic deformations of hyperbolic
3-manifolds. In particular, estimates relating the decay of the norm
of the strain field induced by a deformation to the depth in the
convex core have been absent from prior treatments. The second portion
of the paper develops geometric limit arguments vis a vis the complex
of curves.  The paper concludes with our proof of the convergence of
pseudo-Anosov iteration and double-iteration on quasi-Fuchsian space,
exhibiting explicitly the hyperbolic structure on the pseudo-Anosov
mapping torus $T_\psi$.

\bold{Acknowledgements.}  The authors gratefully acknowledge the
support of the National Science Foundation.  The first author thanks
Guggenheim Foundation and the second thanks the Sloan Foundation for
their support.  We thank MSRI for their hospitality while this work
was being completed, and we thank the referee for many useful
comments.

\section{Deformations}
Let $M$ be a 3-manifold and $g_t$ a one-parameter family of hyperbolic
metrics on $M$ with $D_t$ the covariant derivative for the Riemannian
connection for $g_t$. At time $t=0$ we let $g=g_0$ and $D = D_0$. We
define the time zero derivative, $\eta$, of $g_t$ by the formula
$$\frac{dg_t(v,w)}{dt} \left|_{t=0} \right. = 2g(\eta(v), w).$$
Then $\eta$ is a symmetric tensor of type $(1,1)$. We define the
pointwise norm of $\eta$ at $p$ by choosing an orthonormal basis
$\{e_1, e_2, e_3\}$ for $T_p M$ in the $g$-metric and setting
$$\|\eta\|^2 = \sum_i g(\eta(e_i), \eta(e_i)).$$
Note that this $L^2$-norm bounds the sup norm
%% JFB
from above
so that we have
$$\|\eta(v)\| \leq \|\eta\|\|v\|.$$
If $\eta_t$ is the time $t$ derivative of $g_t$ and $\|\eta_t\| \leq
K$ for all $t \in [0,T]$ then by integrating we see that
$$e^{-2KT}g(v,v) \leq g_T(v,v) \leq e^{2KT}g(v,v).$$
In particular, the identity map on $M$ is a $e^{KT}$-bi-Lipschitz map
from $(M, g)$ to $(M, g_T)$.

We can also use $\eta$ to bound the change in the complex length of
geodesics. Let $\gamma$ be an essential closed curve in $M$ and let
$\cL_\gamma(t) = \ell_\gamma(t) + \imath \theta_\gamma(t)$ be the
complex length of the holonomy of $\gamma$ in the $g_t$-metric.  The
following proposition is a combination of Proposition 4.3 and Lemma
4.6 in \cite{Bromberg:Schwarzian}.  \marginpar{\tiny do I need to say
  more about what complex length is? k}
\begin{prop}
\label{lengthderivative}
Let the harmonic strain field $\eta$ be the time zero derivative of a
family of hyperbolic metrics $M_t = (M, g_t)$. Let $\gamma$ be an
essential simple closed curve in $M$ and $\cL_\gamma(t) =
\ell_\gamma(t) + \imath \theta_\gamma(t)$ its complex length in
$M_t$. Let $\gamma^*$ be the geodesic representative of $\gamma$ in $M_0$.
\begin{enumerate}
\item If the pointwise norms of $\eta$ and $D \eta$ are bounded by $K$ on $\gamma^*$ then
$$|\cL'_\gamma(0)| \leq \sqrt{\frac{2}{3}} K \ell_\gamma(0).$$
\item If $\gamma^*$ has a tubular neighborhood $U$ of radius $R$ then
$$\int_U \|\eta\|^2 + \|D\eta\|^2 \geq \left(\frac{|\cL_\gamma'(t)|}{2\ell_\gamma(t)}\right)^2\left(\frac{\sinh R}{\cosh R}\right) \left(2 + \frac{1}{\cosh^2 R} \right) \area \del U.$$
\end{enumerate}
\end{prop}

When the derivative $\eta$ is a {\em harmonic strain field} there are
a number of formulas that are very useful in controlling the norm of
$\eta$. Before stating these formulas we define harmonic. Given a
family of hyperbolic metrics $(M, g_t)$ around each point we can find
a one-parameter family of $\hthree$-charts $(U, \phi_t)$ for the hyperbolic
structure induced by the $g_t$-metric. These charts can be viewed
as a flow on a neighborhood in $\hthree$. Let $v$ be the vector field
on $U$ that is the pullback of the time zero derivative of this
flow. We then observe $\sym D v = \eta$. This follows from the fact
that for vector fields $u$ and $w$ on $M$ the derivative
$\frac{dg_t(u,w)}{dt}\left|_{t=0} \right.$ is exactly the Lie
derivative of $g(u,w)$ along the vector field $v$.  

The trace of $\sym D v$ is the divergence of $v$ and it measures the
infinitesimal change in volume. The traceless part, $\sym_0 D v$, is
the {\em strain} of $v$ and it measures the infinitesimal change in
the conformal structure. The vector field $v$ is harmonic if
$$D^* D v + 2 v = 0.$$
Here $D^*$ is the formal adjoint of $D$. The factor of $2$ arises from
the fact that the Ricci curvature of a hyperbolic manifold is $-2$,
and the normalization guarantees that infinitesimal isometries are harmonic. We
%% JFB
say
%% see 
that a strain field $\eta$ is harmonic if locally there is a
divergence free and harmonic vector field $v$ with $\eta = \sym D v$.

Finally we note that if $\eta$ is a harmonic strain field then $*D\eta$ is also an harmonic strain field where $*$ is the Hodge star-operator (see Proposition 2.6 in \cite{Hodgson:Kerckhoff:rigidity}). While we are only really interested in controlling the size of $\eta$ we will see throughout the paper that our formulas will also involve $*D\eta$ and we will also control its size along the way.

\section{Infinitesimal inflexibility}
The following formula is our key tool for calculating the $L^2$-norm of  a harmonic strain field. It is Proposition 1.3 of \cite{Hodgson:Kerckhoff:rigidity} along with the calculations on p. 36 of the same paper.

\begin{prop}[Hodgson-Kerkchoff]
\label{boundaryformula}
Let $M$ be a compact manifold with piecewise smooth boundary and $\eta$ a harmonic strain field. Then
$$\int_M \|\eta\|^2 + \|D\eta\|^2 = \int_{\del M} *D\eta \wedge \eta.$$
\end{prop}

The following inequality will allow us to control the boundary term in terms of point-wise bounds on the norms of $\eta$ and $D\eta$.
\begin{lemma}
\label{inequality}
We have $\|\eta\|^2 + \|D\eta\|^2 \geq 2\|*D\eta \wedge \eta\|.$
\end{lemma}

\bold{Proof.} The inequality follows from the fact that $\|\eta - *D\eta\|^2\geq 0$. \qed{inequality}

The following lemma is the first step in showing that the formula from
Proposition~\ref{boundaryformula} holds on some non-compact manifolds
if the strain field is bounded.
\begin{lemma}
\label{exhaust}
Let $M$ be a a complete hyperbolic 3-manifold that is exhausted by compact submanifolds $M_n$ with the area of $\del M_n$ bounded above. If $\eta$ is a harmonic strain field with the pointwise norms $\|\eta\|$ and $\|D\eta\|$ bounded above then the $L^2$ norm of $\eta$ and $D\eta$ is finite.
\end{lemma}

\bold{Proof.} By Proposition~\ref{boundaryformula}
$$\int_{M_n} \|\eta\|^2 + \|D\eta\|^2 = \int_{\del M_n} *D\eta \wedge \eta.$$
Since both the area of $\del M_n$ and the pointwise norms of $\eta$
and $D\eta$ are bounded, Lemma~\ref{inequality} implies that the right
hand side is bounded. This implies that the $L^2$-norm on $M$ is
finite. \qed{exhaust}

Let $P_n$ be a finite $1/n$-net on $\del M$. Define
$$M(t) = \{p \in M | d(p, \del M) \geq t\}$$
and
$$M_n(t) = \{p \in M | d(p, P_n) \geq t\}.$$

\begin{lemma}
\label{nicemanifold}
For all but an isolated set of $t>1/n$, $M_n(t)$ is a manifold with piecewise smooth boundary.
\end{lemma}

\bold{Proof.} If the boundary of $M_n(t)$ is not a manifold with piecewise smooth boundary then there is a geodesic of length $2t$ in $M$ with endpoints in $P_n$. The set of lengths of geodesics in $M$ with endpoints in $P_n$ is a discrete subset of $\reals$ so $M_n(t)$ must be a manifold with piecewise smooth boundary for all but an isolated set of values for $t$. \qed{nicemanifold}

\begin{lemma}
\label{boundaryterm}
Let $M$ be a hyperbolic 3-manifold with piecewise smooth, compact boundary and let $\eta$ be a harmonic strain field on $M$. If $\eta$ and $D\eta$ have finite $L^2$-norm on $M$ then
$$ \int_{M} \|\eta\|^2 + \|D\eta\|^2 = \int_{\del M} *D\eta \wedge \eta.$$
\end{lemma}

\bold{Proof.} Fix a net $P_n$ and a $T>0$ such that the $T$-neighborhood of $P_n$ contains $\del M$ and $M_n(T)$ is a manifold with piecewise smooth boundary. If we apply Proposition~\ref{boundaryformula} to $M \backslash M_n(T)$ and rearrange terms we have
$$\int_M \|\eta\|^2 + \|D\eta\|^2 = \int_{\del M} *D\eta \wedge \eta - \int_{\del M_n(T)} *D\eta \wedge \eta + \int_{M_n(T)} \|\eta\|^2 + \|D\eta\|^2.$$

By Lemma~\ref{nicemanifold} we can choose a sequence of $t_i \rightarrow \infty$ such that $M_n(t_i)$ is a manifold with piecewise smooth boundary. We now apply Proposition~\ref{boundaryformula} again to see that
$$\int_{M_n(T)} \|\eta\|^2 + \|D\eta\|^2 = \int_{\del M_n(T)} *D\eta \wedge \eta - \underset{i \rightarrow \infty}{\lim} \int_{\del M_n(t_i)} *D\eta \wedge \eta.$$
The function
$$f(t) = \int_{\del M_n(t)} (\|\eta\|^2 + \|D\eta\|^2) dA$$
is defined for all but a discrete set of $t$ and therefore
$$\int_{M_n(T)} \|\eta\|^2 + \|D\eta\|^2 = \int_T^\infty f(t) dt.$$
Since the $L^2$-norm of $\eta$ and $D\eta$ is finite on $M_n(T)$ we have
$$\underset{t \rightarrow \infty}{\lim}f(t) = 0$$
and in particular $f(t_i) \rightarrow 0$.
Lemma~\ref{inequality} then implies that
$$f(t_i) \geq 2 \left| \int_{\del M_n(t_i)} *D\eta \wedge \eta \right|.$$
Therefore
$$\underset{i \rightarrow \infty}{\lim} \int_{\del M_n(t_i)} *D\eta \wedge \eta = 0$$
and
$$\int_{M_n(T)} \|\eta\|^2 + \|D\eta\|^2 = \int_{\del M_n(T)} *D\eta \wedge \eta.$$
Combining this last equality with the first equality in the proof gives us the lemma.
\qed{boundaryterm}

\bold{Remark.} Sullivan's rigidity theorem, which guarantees that a
quasiconformal deformation of a finitely generated Kleinian group
$\Gamma$ with support in the limit set is trivial, played a central
role in Thurston's original proof of the existence of hyperbolic
structures on fibered 3-manifolds.  Thurston \cite{Thurston:book:GTTM}
and Bonahon \cite{Bonahon:tame} subsequently observed that Sullivan
rigidity \cite{Sullivan:linefield} follows somewhat more directly if
one assumes the tameness of $M = \half^3/\Gamma$, namely that $M$ is
homeomorphic to the interior of a compact 3-manifold
(cf. \cite[\S3]{McMullen:book:RTM}).

Lemmas
\ref{exhaust} and \ref{boundaryterm} give another perspective on
Sullivan's result.  
In particular, 
% if $M = \hthree/\Gamma$ is a
% complete hyperbolic 3-manifold then 
%% JFB
%% limit
any $\Gamma$-invariant Beltrami differential extends continuously via
an averaging process to a harmonic strain field $\eta$ on $M$ with the
pointwise norms of $\eta$ and $D\eta$ uniformly bounded.  If $M$ is
tame, then the limit set of $\Gamma$ has measure zero or is all of
$\chat$, by Canary's \cite{Canary:ends} result that tameness implies
Ahlfors' measure conjecture.  In the former case, any Beltrami
differential supported on the limit set is trivial. In the later case
tameness also implies that $M$ is exhausted by submanifolds whose
boundary has uniformly bounded area and Lemma \ref{exhaust} implies
that the $L^2$-norms of $\eta$ and $D\eta$ are finite on $M$. Since
$M$ has no boundary, Lemma \ref{boundaryterm} implies that $\eta =
D\eta = 0$ and the initial Beltrami differential must be trivial.

\medskip

The following theorem is the key analytic estimate that underlies all of our inflexibility theorems. It should be compared with Theorem 2.15 in \cite{McMullen:book:RTM}.
\marginpar{\tiny Give this a name? Infinitesimal inflexibility? The reference should be to CTM's theorem of the same name. k}
\begin{theorem}
\label{expdecay}
Let $M$ be a hyperbolic 3-manifold with compact boundary and let
$\eta$ be a harmonic strain field on $M$.  Assume that the $L^2$-norms
$\|\eta\|^2$ and $\|D\eta\|^2$ are finite.
Then
$$\int_{M(t)} \|\eta\|^2 + \|D\eta\|^2 \leq e^{-2t}\int_M \|\eta\|^2 + \|D\eta\|^2.$$
\end{theorem}

\bold{Proof.}  We will show that
\begin{equation}
\label{theineq}
\int_{M_n(t)} \|\eta\|^2 + \|D\eta\|^2 \leq e^{-2(t - 2/n)} \int_{M_n(2/n)} \|\eta\|^2 + \|D\eta\|^2.
\end{equation}
Taking the limit of this inequality as $n \rightarrow \infty$ will imply the theorem.

Let
$$f(t) = \int_{M_n(t)} \|\eta\|^2 + \|D\eta\|^2.$$
By Lemma~\ref{nicemanifold} we can write
$$f(T) = \int_T^\infty \int_{\del M_n(t)} \left( \|\eta\|^2 + \|D\eta\|^2\right) dA dt$$
for $T > 1/n$. Therefore
\begin{eqnarray}
-f'(t) & = & \int_{\del M_n(t)} \left(\|\eta\|^2  + \|D\eta\|^2\right) dA \nonumber \\
& \geq & 2 \int_{\del M_n(t)} *D\eta \wedge \eta \nonumber \\
& \geq & 2f(t). \nonumber
\end{eqnarray}
Integrating both sides of this inequality from $2/n$ to infinity implies \eqref{theineq}. \qed{expdecay}

To go from $L^2$-bounds on $\eta$ to pointwise bounds we use the following mean value theorem of Hodgson and Kerckhoff. A proof can be found in \cite{Bromberg:Schwarzian}.

\begin{theorem}
\label{meanvalue}
Let $\eta$ be a harmonic strain field on a ball $B$ of radius $R< \pi/2$ centered at a point $p$. Then
$$\|\eta(p)\| \leq \frac{3\sqrt{2\vol{B}}}{4\pi f(R)} \sqrt{\int_B \|\eta\|^2}$$
where $f(R) = \cosh(R) \sin(\sqrt{2}R) - \sqrt{2}\sinh(R) \cos(\sqrt{2}R)$.
\end{theorem}

We also recall the Margulis {\em thick-thin decomposition} for
hyperbolic surfaces and 3-manifolds. If $M$ is a Riemannian manifold
{\em injectivity radius} $\inj_M\colon M \to \reals^+$ measures the
radius of the maximal embedded metric ball at each point.  The {\em
  $\epsilon$-thin part} of $M$, denoted $M^{< \epsilon}$, is the set
of points $x$ in $M$ for which $\inj_M(x) < \epsilon$.  Likewise the
{\em $\epsilon$-thick part} $M^{\ge \epsilon}$ of $M$ is the set of
$x$ for which $\inj_M(x) \ge \epsilon$.
\begin{lemma}[Margulis] {\sc (Thick-Thin Decomposition)}
\label{Margulis}
There exists $\varepsilon_2 >0$ so that if $X$ is a hyperbolic surface
and $\epsilon \leq \varepsilon_2$ then every component of
$X^{<\epsilon}$ is either the open metric $R$-neighborhood of a simple closed
geodesic, $R>0$, or an open horosphere modulo a
discrete parabolic $\zed$ action.

There exists $\varepsilon_3 >0$ so that if $M$ is a complete
hyperbolic 3-manifold and $\epsilon \leq \varepsilon_3$ then every
component of $M^{< \epsilon}$ is either the open metric tubular
$R$-neighborhood $\tube_\epsilon(\gamma)$ of a simple closed geodesic
$\gamma$ in $M$, or an open horoball modulo a discrete parabolic $\zed
$ or $\zed \oplus \zed$ action.
\end{lemma}
The tube $\tube_\epsilon(\gamma)$ is called a {\em Margulis tube}, and
the horoball quotients are called {\em rank-1} or {\em rank-2} cusps
depending on whether the action is by a $\zed $ or $\zed \oplus \zed$
parabolic subgroup of $\PSL_2(\cx)$.  We employ the notation
$\tube_{\varepsilon_3}(\gamma) = \tube(\gamma)$.

% It is useful to let $M^{\ge \epsilon}$ refer to the locus of points in
% $M$ where the injectivity radius is at least $\epsilon$.  For complete
% hyperbolic 3-manifolds, the Margulis {\em thick-thin decomposition}
% for hyperbolic $3$-manifolds guarantees the existence of a universal
% $\epsilon_3>0$ (depending only on the dimension) for which each
% portion of the $\epsilon_3$-thin part, where the injectivity radius is
% less than $\epsilon_3$, has a the standard form of a solid torus
% neighborhood $\tube$ of a short geodesic $\gamma^*$ (a {\em
%   Margulis-tube}), a horoball modulo a parabolic $\zed$-action (a {\em
%   rank-one cusp}) or a horoball modulo a parabolic $\zed \dirsum
% \zed$-action (a {\em rank-two cusp}).

We now apply Theorems~\ref{expdecay} and~\ref{meanvalue} to obtain pointwise bounds on $\eta$.

\begin{theorem}
\label{pointdecay}
Let $M$ be a complete hyperbolic 3-manifold with compact boundary and let $\eta$ be a harmonic strain field on $M$. Then
$$\|\eta(p)\| \leq A(\epsilon)e^{-d(p, \del M)} \sqrt{\int_M \|\eta\|^2 + \|D\eta\|^2}$$
where $p \in M^{\geq \epsilon}$ and
$$A(\epsilon) = \frac{3e^{\epsilon}\sqrt{2 \vol(B)}}{4\pi f(\epsilon)}$$
with the function $f$ defined in Theorem~\ref{meanvalue}.
\end{theorem}

\bold{Proof.} Let $B$ be the ball of radius $\epsilon$ centered at
$p$. Then $B$ lies in $M(d(p, \del M) - \epsilon)$. By Theorem
\ref{expdecay}
\begin{eqnarray*}
\int_B \|\eta\|^2 & \leq & \int_B \|\eta\|^2 + \|D\eta\|^2 \\
& \leq & \int_{M(d(p, \del M) - \epsilon)} \|\eta\|^2 + \|D\eta\|^2 \\
& \leq & e^{-2(d(p,\del M) - \epsilon)} \int_M \|\eta\|^2 + \|D\eta\|^2
\end{eqnarray*}
We then apply Theorem~\ref{meanvalue} to finish the proof. \qed{pointdecay}

We can also control the derivative of the length of a closed geodesic.

\begin{theorem}
\label{lengthdecay}
Let the harmonic strain field $\eta$ be the time zero derivative of a
family of hyperbolic metrics $M_t = (M, g_t)$ where $M$ is a
3-manifold with compact boundary. Let $\gamma$ be an essential simple
closed curve in $M$ and $\cL_\gamma(t) = \ell_\gamma(t) + \imath
\theta_\gamma(t)$ its complex length in $M_t$. Let $\gamma^*$ be the
geodesic representative of $\gamma$ in $M_0$. 
\begin{enumerate}
\item If $\gamma^*$ is contained in $M^{\geq \epsilon}_0$ then
$$|\cL'_\gamma(0)| \leq A(\epsilon)e^{-d(\gamma^*, \del M)} \ell_\gamma(0) \sqrt{\frac{2}{3}\int_M \|\eta\|^2 + \|D\eta\|^2}$$
where $A(\epsilon)$ is the function given in Theorem~\ref{pointdecay}.
\item If $\gamma^*$ has a tubular neighborhood $U$ of radius $R$ then
$$|\cL'_\gamma(0)| \leq C(R) e^{-d(U, \del M)}\ell_\gamma(0) \sqrt{\frac{ \int_M \|\eta\|^2 + \|D \eta\|^2} {\area(\del U)}}$$
where $1/C(R) =2 \tanh R \left(2 + \frac{1}{\cosh^2 R}\right)$.
\end{enumerate}
\end{theorem}

\bold{Proof of (1).} Applying Theorem~\ref{pointdecay} we see that on $\gamma^*$ the pointwise norm of $\eta$ and $D\eta$ is bounded by $A(\epsilon) e^{-d(\gamma^*, \del M)} \sqrt{\|\eta\|^2 + \|D\eta\|^2}$. We then apply (1) of Proposition~\ref{lengthderivative} to finish the proof.

\bold{Proof of (2).} By Theorem~\ref{expdecay}
$$\int_U \|\eta\|^2 + \|D\eta\|^2 \leq e^{-d(U, \del M)} \int_M \|\eta\|^2 + \|D\eta\|^2.$$
In this case, (2) of Proposition~\ref{lengthderivative} finishes the
proof. \qed{lengthdecay}

\section{Inflexibility}
\label{section:inflexibility}
There are two types of deformations of hyperbolic 3-manifolds that can be
studied with our methods: quasiconformal deformations, namely, quasiconformal
conjugacies of their uniformizing Kleinian groups, and cone-manifold
deformations, deformations of a singular {\em cone-manifold} structure
wherein the cone-angle at the cone-locus varies. In this paper we will restrict to quasiconformal deformations deformations but the two general global inflexibility theorems we prove in this section can also be applied to the study of cone-manifolds. We will carry this out in a subsequent paper.

\begin{theorem}
\label{maininflex}
Let $g_t$ be a one-parameter family of hyperbolic metrics on a
3-manifold $M$ with $t \in [a,b]$.  Let $\eta_t$ be the time $t$
derivative of the metrics $g_t$ and let $N_t$ be a family of 
%% JFB
3-dimensional
submanifolds of $M$ such that $\eta_t$ is a harmonic strain field on
$N_t$. Also assume that
$$\sqrt{\int_{N_t} \|\eta_t\|^2 + \|D_t\eta_t\|^2} \leq K$$
for some $K>0$.
Let $p$ be a point in $M$ such that for all $t\in [a,b]$, $p$ is in $M_t^{\geq \epsilon}$ and
$$d_{M_t}(p, M\backslash N_t) \geq	d$$
where $d>\epsilon$.
Then
$$\log \bilip(\Phi_t, p) \leq (t-a)KA(\epsilon)e^{-d}$$
where $\Phi_t$ is the identity map from $M_a$ to $M_t$ and $A(\epsilon)$ is the function from Theorem~\ref{pointdecay}.
\end{theorem}
\marginpar{\tiny Need to decide where to define $M^{\geq \epsilon}$, etc... k}

\bold{Proof.} Since $d>\epsilon$ the $\epsilon$-neighborhood of $p$ is
contained in $N_t$ and is at least distance $d-\epsilon$ from $\del
N_t$. An application of Theorem~\ref{pointdecay} gives us
$$\|\eta_t(p)\| \leq KA(\epsilon)e^{-d}.$$
Integrating we get
$$\log \bilip(\Phi_t, p) \leq (t-a)KA(\epsilon)e^{-d}$$
as desired. \qed{maininflex}

Though the previous result gives no control over the bi-Lipschitz
constant of the map $\Phi$ in the thin part, we may instead
demonstrate exponential decay of the change in length of short curves
in Margulis thin parts, which controls the geometry of the thin part
itself.  Here, the decay is measured in terms of the distance of the
corresponding Margulis tube from the boundary. For completeness, we
also bound the change in length of curves that are not short.

\begin{theorem}
\label{mainlengthinflex}
Let $g_t$ be a one-parameter family of hyperbolic metrics on a
3-manifold $M$ with $t \in [a,b]$.  Let $\eta_t$ be the time $t$
derivative of the metrics $g_t$ and let $N_t$ be a family of
submanifolds of $M$ such that $\eta_t$ is a harmonic strain field on
$N_t$. Also assume that
$$\sqrt{\int_{N_t} \|\eta_t\|^2 + \|D_t\eta_t\|^2} \leq K$$
for some $K>0$. Let $\gamma_t$ be the geodesic representative on $(M,
g_t)$ of a closed curve $\gamma$ and let $\ell_\gamma(t)$ be the
length of $\gamma$.
\begin{enumerate}
\item Assume that $\gamma_t$ is in $M_t^{\geq \epsilon}$ for all $t\in
  [a,b]$, and that
$$d_{M_t}(\gamma_t, M\backslash N_t) \geq	d.$$
Then
$$\left|\log \frac{\ell_\gamma(b)}{\ell_\gamma(a)}\right| \leq \sqrt{2/3}A(\epsilon)(b-a)Ke^{-d}.$$
\item Assume $\gamma_t$ has a tubular neighborhood $U_t$ of radius
  at least $R$ and the area of $\del U_t$ is  at least $B$. Also assume that
$$d_{M_t}(U_t, M\backslash N_t) \geq d$$
for all $t\in [a,b]$.
Then
$$\left|\log \frac{\ell_\gamma(b)}{\ell_\gamma(a)}\right| \leq \frac{C(R)(b-a)Ke^{-d}}{\sqrt{B}}$$
where $C(R)$ is the function from Theorem~\ref{lengthdecay}.
\end{enumerate}
\end{theorem}

\bold{Proof.} Both inequalities are obtained by integrating the
estimates of Theorem~\ref{lengthdecay}. \qed{mainlengthinflex}

\bold{Remark.} Although in the above theorem we only control the real
lengths of closed geodesics it is straightforward to control their
complex lengths.  In particular if $\cL_\gamma(t)$ is the complex
length of $\gamma$ in $(M,g_t)$ then we can view $\imath
\cL_\gamma(t)$ as a point in the upper half space model of
$\htwo$. Then the quantities on the right hand side of the
inequalities bound the hyperbolic distance between $\imath
\cL_\gamma(a)$ and $\imath \cL_\gamma(b)$. Note that this hyperbolic
distance is an upper bound on the $\log$ of the ratio of real lengths
so such a hyperbolic distance bound implies the inequalities in
Theorem~\ref{mainlengthinflex}.

\section{Quasiconformal deformations}
We now apply the results of the previous section to quasiconformal
deformations. We begin reviewing some standard definitions.

Let $M$ be a complete, orientable, hyperbolic 3-manifold. Its
universal cover is naturally identified with $\hthree$ and $M$ may be
recovered as the quotient $M = \hthree / \Gamma$ of $\hthree$ by a
Kleinian group $\Gamma$, namely, a discrete subgroup of
$\Isom^+(\hthree)$. The natural action of $\Gamma$ on $\chat$ by
M\"obius transformations partitions $\chat$ into its {\em domain of
  discontinuity}, $\Omega$, the largest subset of $\chat$ where
$\Gamma$ acts properly discontinuously, and its {\em limit set}
$\Lambda$. Then the {\em Kleinian manifold} quotient $(\hthree \cup
\Omega)/\Gamma $ is a 3-manifold with {\em conformal boundary}
$\Omega/\Gamma$.

A {\em $K$-quasiconformal deformation} of a complete, orientable,
hyperbolic 3-manifolds $M_0$ is a is a map $\Psi: M_0 \to
M_1$ to a complete hyperbolic 3-manifold $M_1$ such that the lift
$\tilde{\Psi} : \hthree \to \hthree$ to the universal
covers extends continuously to a $K$-quasiconformal map of $\chat$. If
$\Psi$ is a $K$-quasiconformal deformation then it will extend to a
$K$-quasiconformal map between the conformal boundaries of $M_0$ and
$M_1$.
\marginpar{\tiny do you want $\Psi$ to be a diffeo? -J}

The following result is due to Reimann \cite{Reimann:visual} using
work of Ahlfors \cite{Ahlfors:visual} and Thurston
\cite{Thurston:book:GTTM}. For a self contained exposition see
\cite{McMullen:book:RTM}. It is an essential tool for the work that
follows
\begin{theorem}[Reimann]
\label{qcflow}
Let $\Psi: M_0 \to M_1$ be a $K$-quasiconformal
deformation of the complete hyperbolic 3-manifold $M_0$. Then
there exists a one-parameter family, $M_t = (M, g_t)$, $t \in [0,1]$,
of hyperbolic 3-manifolds with time $t$ derivative $\eta_t$ such that
the following holds:
\begin{enumerate}
\item The $\eta_t$ are harmonic strain fields and $\|\eta_t\|_\infty,
  \|D_t \eta_t\|_\infty \leq 3k$ where $k = \frac{1}{2} \log K$;

\item Let $\Phi_t : M_0 \to M_t$ be the identity map on
  $M$. Then $\Phi_t$ is $K^{\frac{3}{2}}$-bi-Lipschitz and $\Phi_1$ is
  homotopic to $\Psi$.
\end{enumerate}
\end{theorem}

The convex
cores $C(M_t)$ of the one-parameter family $M_t$ will play the role
of $N_t$ when we apply Theorems~\ref{maininflex} and~\ref{mainlengthinflex} to $M_t$.

\begin{lemma}
\label{finiteconvex}
Let $M$ be a complete hyperbolic 3-manifold such that $\pi_1(M)$ is
finitely generated and assume that $M$ has no rank one cusps. Let
$\eta$ be a harmonic strain field on $M$ such that the norms of $\eta$
and $D \eta$ are pointwise bounded by $k$. Then
$$\int_{C(M)} \|\eta\|^2 + \|D \eta\|^2 \leq \area(\del C(M)) k^2.$$
\end{lemma}

\bold{Proof.} We first replace the convex core with its
$\epsilon$-neighborhood, $C_\epsilon(M)$. While the boundary of the
convex core may not be smooth, the boundary of $C_\epsilon(M)$ will be
$C^1$. We also note that $\area(\del C_\epsilon(M)) \rightarrow
\area(\del C(M))$ as $\epsilon \rightarrow 0$.  \marginpar{\tiny The
  fact about being $C^1$ is due to Bowditch or at least there is a
  nice proof due to him. I think the proof is at the start of the
  original Epstein-Marden paper. I don't know a reference for the
  areas limiting. k}

Since $\pi_1(M)$ is finitely generated the $M$ are both topologically
and geometrically tame (\cite{Agol:tame, Calegari:Gabai:tame}). In
particular the convex cores $C(M)$ will be exhausted by submanifolds
whose boundary has uniformly bounded area. Since the norms of $\eta$
and $D \eta$ are uniformly bounded we can apply Proposition
\ref{boundaryformula} and Lemma~\ref{inequality} to see that the
$L^2$-norms of $\eta$ and $D \eta$ are uniformly bounded on these
submanifolds which implies that the $L^2$-norms of $\eta$ and $D \eta$
are finite on $C_\epsilon(M)$.

Applying Lemma~\ref{boundaryterm} to $C_\epsilon(M)$ and taking a limit as $\epsilon \rightarrow 0$ gives us the lemma. \qed{finiteconvex}

To make sure that objects deep in the convex core of $C(M_0)$ stay
deep in the convex core of $C(M_1)$ we will use the fact that
bi-Lipschitz maps of $\half^3$ take convex subsets of $\half^3$ to
quasi-convex sets, a general feature of quasi-isometries between
$\delta$-hyperbolic spaces. While this section only applies this
observation for hyperbolic space we will later make use of this more
general version in the setting of manifolds with pinched negative
curvature, so we give more general form.  Sometimes known as the Morse
Lemma, Theorem 1.7 in Chapter of III.H of \cite{Bridson:Haefliger:npc}
is one reference.  \marginpar{\tiny Should we say more/less about the
  history of this fact? I don't really care if we call it the Morse
  lemma. k}

\begin{theorem}\label{morselemma}
  Given constants $L>1$ and $\epsilon \in (0,1)$ there exists a $K>0$
  such that the following holds.  Let $X_0$ and $X_1$ be complete,
  simply connected Riemannian manifolds with sectional curvatures
  lying in $(-1 -
  \epsilon, -1 + \epsilon)$, and let $\Phi: X_0 \to X_1$ be an
  $L$-bi-Lipschitz diffeomorphism. Then the $\Phi$-image of a convex
  set in $X_0$ is $K$-quasi-convex in $X_1$.
\end{theorem}

An example of a convex set is a geodesic -- its image under
a bi-Lipschitz map is an example of a quasi-geodesic. A more common
way to state the above theorem is that in a space with pinched
negative curvature,  a
quasi-geodesic is a bounded Hausdorff distance from a geodesic. In
fact this is how the result is stated in \cite{Bridson:Haefliger:npc}
but it is not hard to see that this implies the above theorem.

On application of the above theorem is the following proposition.
\begin{prop}\label{coresclose}
  Given $B>1$ and $\epsilon \in (0,1)$ there exists $d>0$ such that
  the following holds. Let $g_0$ and $g_1$ be complete Riemannian
  metrics on a manifold $M$ with sectional curvatures in
  $(-1-\epsilon,-1 + \epsilon)$ and let $\phi: (M, g_0)
  \to (M, g_1)$ be $B$-bi-Lipschitz. Then then Hausdorff
  distance between $C(M,g_1)$ and $\phi(C(M,g_0))$ is less than $d$.
\end{prop}
\marginpar{\tiny Say something about the definition of the convex core of a negatively curved manifold? k}

\bold{Proof.} For hyperbolic manifolds this is Proposition 2.16 in
\cite{McMullen:book:RTM}. It follows from Theorem~\ref{morselemma} and
the fact that every point in the the convex hull of a set is a uniform
distance from a geodesic with endpoints in the set. Using work of
Anderson \cite{Anderson:dirichlet}, Bowditch \cite{Bowditch:hulls}
proved this last fact for manifolds with pinched negative curvature
where the uniformity constants depend on the pinching constants. Using
Bowditch's work, McMullen's proof extends to the setting we have
here. \qed{coresclose}

The following is Corollary 2.17 in \cite{McMullen:book:RTM}. The proof is a straightforward application of  Proposition~\ref{coresclose}.
\begin{lemma}
\label{staydeep}
Let $\Phi: M_0 \to M_1$ be an $L$-bi-Lipschitz diffeomorphism between complete hyperbolic 3-manifolds. Then there exist a constant $d$ such that
$$d(\Phi(p), M_1 \backslash \core(M_1)) \geq \frac{d(p, M_0 \backslash \core(M_0))}{L} - d.$$
\end{lemma}
\marginpar{\tiny Curt uses $d$ so I used $d$ but is that the best choice? k}

We are now ready to prove our first inflexibility theorem for quasiconformal deformations.
\begin{theorem}
\label{thickinflex}
Let $M_0$ and $M_1$ be complete hyperbolic structures on a 3-manifold
$M$ such that $M_1$ is a $K$-quasiconformal deformation of $M_0$,
$\pi_1(M)$ is finitely generated, and $M_0$ has no rank one cusps.
Then there is a bi-Lipschitz diffeomorphism
$$\Phi \colon M_0 \to M_1$$
whose pointwise bi-Lipschitz constant satisfies
$$\log \bilip(\Phi, p) \le C_1e^{-C_2 d(p, M_0 \backslash \core(M_0))}$$
where $p$ is in $M_0^{\ge \epsilon}$ and $C_1$ and $C_2$ depend only on $K$,
$\epsilon$, and $\area(\del\core(M_0))$.
\end{theorem}

\bold{Proof.} Let $M_t = (M, g_t)$ be the one-parameter family of
hyperbolic manifolds given by 
Theorem~\ref{qcflow} 
with $\eta_t$ the
derivative of the metrics and $$
\Phi_t \colon M_0 \to M_t
$$ the
given maps. By Lemma~\ref{finiteconvex}  we have
$$\int_{C(M_t)} \|\eta_t\|^2 + \|D_t \eta_t\|^2 \leq \area(\del C(M_t))9k^2.$$
Lemma~\ref{staydeep} guarantees
$$d(\Phi_t(p), M_t \backslash C(M_t)) \geq \frac{d(p, M_0 \backslash C(M_0))}{K^{\frac{3}{2}}} - d.$$
Since 
by Theorem~\ref{qcflow} 
the $\Phi_t$ are $K^{\frac{3}{2}}$-bi-Lipschitz we have $p \in
M_t^{\ge \epsilon'}$ for all $t$ where $\epsilon' =
\epsilon/K^{\frac{3}{2}}$. The result then follows from Theorem
\ref{maininflex} with $\Phi = \Phi_1$ the desired
map. \qed{thickinflex}

For points in the thin part, the above theorem fails to yield good
estimates, but this is not surprising.  Indeed, one can construct
examples of harmonic strain fields on Margulis tubes where the
pointwise $L^2$-norm is roughly constant and does not decay with depth
into the tube. Rather, one expects the pointwise norm of the strain at
a point in a Margulis tube to depend on the depth of the boundary of
the tube. Rather than pursue such a line of argument, we will bound
the change in length of short geodesics where, again, the bounds will
depend on the depth of the boundary of the Margulis tube not the short
geodesic. Such a bound is the natural thing to expect and suffices for
applications.  \marginpar{\tiny Do we want any of this commentary in
  the intro? It does emphasize that our theorem is somehow optimal in
  that we can only hope to have exponential decay of the bilip
  constant for points in the thick part. k}

For completeness we also give bounds on the change in length of curves
that have bounded length but are not necessarily short.  We must first
show that an essential curve whose geodesic representative lies deep in
the convex core of $M_0$ also has geodesic representative in $M_1$
deep in the convex core.

\begin{prop}\label{tubestaydeep}
  Let $M_0=(M, g_0)$ and $M_1 = (M, g_1)$ be hyperbolic 3-manifolds
  that are $L$-bi-Lipschitz diffeomorphic.  Let $\epsilon$ a positive
  constant such that $L\epsilon < \varepsilon_3$. Then there exists a
  constant $d = d(L, \epsilon)$ such that the following holds. Let
  $\gamma$ be an essential closed curve in $M$ and $\gamma_0$ and
  $\gamma_1$ its geodesic representatives in $M_0$ and $M_1$,
  respectively.
\begin{enumerate}
\item We have $$d(\gamma_1, M_1 - C(M_1)) \geq \frac{d(\gamma_0, M_0 -
    C(M_0))}{L} - d,$$ and,

\item if $\ell_{M_0}(\gamma_0) \leq 2 \epsilon/L$ then
$$d(\tube^1_\epsilon(\gamma), M_1 - C(M_1)) \geq
\frac{d(\tube^0_\epsilon(\gamma), M_0 - C(M_0))}{L} - d$$
where $\tube^1\epsilon(\gamma)$ and $\tube^0\epsilon(\gamma)$ denote
the Margulis tubes about $\gamma$ in $M_1$ and $M_0$.
\end{enumerate}
\end{prop}

\bold{Proof.} Let $\Phi: M_0 \to M_1$ be the
$L$-bi-Lipschitz diffeomorphism. Let $q$ be a point on $\gamma_1$ with
$$d(q, M_1 - C(M_1)) = d(\gamma_1, M_1 - C(M_1)).$$
By Theorem~\ref{morselemma}, the Hausdorff distance between
$\Phi(\gamma_0)$ and $\gamma_1$ is bounded by $K$ where $K$ only
depends on $L$ so there exists a $q' \in \Phi(\gamma_0)$ with $d(q,q')
\leq K$. Let $p = \Phi^{-1}(q')$. Then
$$d(p, M_0 - C(M_0)) \geq d(\gamma_0, M_0 - C(M_0)).$$
An application of Lemma~\ref{staydeep} to $p$ gives us (1).

The proof of (2) is similar with one change. Again let $q$ be a point
on $\del \tube^1_\epsilon(\gamma)$ such that
$$d(q,M_1 - C(M_1)) = d(\tube^1_\epsilon(\gamma), M_1 - C(M_1)).$$
The collar $\tube^1_{L\epsilon}(\gamma) - \tube^1_{\epsilon/L}(\gamma)$ will
contain $\del \Phi(\tube^0_\epsilon(\gamma))$ and the inclusion will be a
homotopy equivalence since $\Phi(\tube^0_\epsilon(\gamma))$ is not
contained in the collar. By \cite{Brooks:Matelski:collars} the width
of the collar is bounded above by some $W$ depending only on
$\epsilon$ and $L$. Therefore there exists a $q' \in \del
\Phi(\tube^0_\epsilon(\gamma))$ such that $d(q,q') \leq W$. The rest of
the proof is the same as in (1).  \qed{tubestaydeep}

We can now control the length of geodesics under quasiconformal deformations.
\begin{theorem}
\label{qclengthinflex}
Let $M_1= (M, g_1)$ be a $K$-quasiconformal deformation of the
hyperbolic 3-manifold $M_0 = (M, g_0)$ with finitely generated
fundamental group and no rank-one cusps. Let $\gamma$ be an essential
simple closed curve in $M$ and $\gamma_0$ and $\gamma_1$ its geodesic
representatives in $M_0$ and $M_1$ respectively. Choose $\epsilon>0$
such that $\epsilon K^{\frac{3}{2}} < \varepsilon_3$, and let $L> 2\epsilon>0$. Then there exists constants $C_1$
and $C_2$ depending on $K$, $\epsilon$, $L$ and $\area(\del C(M_0))$
such that the following holds.
\begin{enumerate}
\item If $2\epsilon \leq \ell(\gamma_0) \leq L$ then
$$\left| \log \frac{\ell(\gamma_1)}{\ell(\gamma_0)}\right| \leq C_1 e^{-C_2 d(\gamma_0, M_0 - C(M_0))}.$$

\item If $\ell(\gamma_0) \leq 2\epsilon$ then
$$\left| \log \frac{\ell(\gamma_1)}{\ell(\gamma_0)} \right| \leq C_1 e^{-C_2 d(\tube_\epsilon^0(\gamma), M_0 - C(M_0))}.$$
\end{enumerate}
\end{theorem}

\bold{Proof.} As with the proof Theorem~\ref{thickinflex} we now only
need to put together the pieces. We will use
Theorem~\ref{mainlengthinflex}, our generic inflexibility theorem for
lengths of curves. To apply this result we use the family of
deformations given by Theorem~\ref{qcflow} where the bound on the
$L^2$-norms of the strain fields inside the convex core comes from
Lemma~\ref{finiteconvex}. Finally, Proposition~\ref{tubestaydeep}
guarantees that geodesics and tubes that are deep in the convex core
stay deep in the convex core. The theorem then follows from an
application of Theorem~\ref{thickinflex}. \qed{qclengthinflex}

\bold{Remark.} It is easy to see that both Theorems~\ref{thickinflex}
and~\ref{qclengthinflex} hold for geometrically finite hyperbolic
manifolds with rank-one cusps. To see this let $M^{< \delta}_c$ be set
of points in the rank one cusps of $M$ that have injectivity radius
less than $\delta$. If $M$ is geometrically finite then $C^\delta_c(M)
= C(M) \backslash M^{< \delta}_c$ will be compact and
Theorems~\ref{thickinflex} and~\ref{qclengthinflex} will hold if we
replace $C(M)$ with $C^\delta_c(M)$. We also note that $\area(\del
C^\delta_c(M)) \rightarrow \area(C(M))$ as $\delta \rightarrow 0$ and
for all $p \in C(M)$ there exists a $\delta_p$ such that if $\delta <
\delta_p$ then 
$$d(p, M\backslash C^\delta_c(M)) = d(p, M \backslash C(M)).$$
 Therefore if we let $\delta \rightarrow 0$ we recover
Theorems~\ref{thickinflex} and~\ref{qclengthinflex} as stated above.
\marginpar{\tiny The only place in the run-up of the proofs
  of~\ref{thickinlex} and ~\ref{qclengthinflex} where we use that
  there aren't rank one cusps is~\ref{finiteconvex}. Should this be
  mentioned or is it implicit?  Also, should this be about the generic
  $M$ or should I say $M_0$, etc... k}

In fact, the above argument applies whenever $C^\delta_c(M)$ is a
manifold with compact boundary, as is the case when either the
intersection of each rank one cusp with the convex core has finite
volume or the entire rank one cusp is contained in the convex core.

We expect both theorems should hold for any hyperbolic 3-manifold with
finitely generated fundamental group.  \marginpar{\tiny Should we see anything about an improved
  inflexibility theorem when there are cusps? Or will we do this later
  when we talk about why the double limit theorem only works for
  closed surfaces? k}

\section{Schwarzian derivatives}
\label{section:schwarzian}
The conformal boundary of a hyperbolic 3-manifold also has a
projective structure. In this section we will obtain bounds on how
this projective boundary changes during a quasiconformal
deformation. We begin with some background on
projective structures. One reference for this material is
\cite{Dumas:handbook}.

A {\em complex projective structure} on a surface $S$ can be defined
in two equivalent ways.  First, a complex
projective structure is an atlas of charts to $\chat$ whose 
transition functions are restrictions of M\"obius transformations. 
Second, a projective structures is a {\em developing pair} $(D, \rho)$
where $D \colon \widetilde{S} \to \chat$ is
a local homeomorphism and $\rho$ is representation of $\pi_1(S)$ in
$\PSL_2\cx$ for which
that $D \circ g(x) = \rho(g) \circ D(x)$ for all $g \in \pi_1(S)$ and
$x \in \widetilde{S}$. The map $D$ is {\em developing map} and $\rho$ is
the {\em holonomy representation}. 
An atlas determines a developing pair and a developing
pair determines an
atlas. 

A projective structure determines a conformal structure on $S$ but
distinct projective structures may have the same underlying conformal
structure. If $X$ is a conformal structure on $S$ then we let $P(X)$
denote the space of projective structures on $S$ with conformal
structure $X$.

Note that the charts that define a conformal structure on the boundary
%% JFB
at infinity
of a hyperbolic 3-manifold also define a projective structure. We
refer to this projective structure as the {\em projective boundary} of the
manifold. We will be interested in controlling how the projective
boundary changes under a deformation fixing the conformal
boundary.

The difference between two projective structures $\Sigma_0$ and
$\Sigma_1$ in $P(X)$ is measured by a quadratic differential $\Phi$
determined via the {\em Schwarzian derivative}. If $f$ is the conformal map between $\Sigma_0$ and $\Sigma_1$ then the Schwarzian derivative of $f$ is the quadratic differential
$$\Phi = \left[\left(\frac{f_{zz}}{f_z}\right)_z - \frac{1}{2}\left(\frac{f_{zz}}{f_z}\right)^2 \right] dz^2$$
where the derivatives are taken in projective charts for $\Sigma_0$ and $\Sigma_1$. We can then define
$d(\Sigma_0, \Sigma_1) = \|\Phi\|_\infty$ where $\|\Phi\|_\infty$ is
the sup-norm taken with respect to the hyperbolic metric on $X$. 

There is also an infinitesimal version of the Schwarzian. If
$\Sigma_t$ is a smooth path in $P(X)$ from $\Sigma_0$ to $\Sigma_1$ then the Schwarzian's from $\Sigma_0$ to $\Sigma_t$ determine a smooth path of quadratic differentials.
The time $t$ derivative $\Phi_t$ of this path is also a quadratic
differential. The following inequality will be useful for bounding
$d(\Sigma_0, \Sigma_1)$:
$$\|\Phi\|_\infty \leq \int_0^1 \|\Phi_t\|_\infty dt.$$

For each hyperbolic structure $X$ there is a unique Fuchsian
projective structure $\Sigma_F$ in $P(X)$. For an arbitrary $\Sigma
\in P(X)$ we define $\|\Sigma\|_F = d(\Sigma,
\Sigma_F)$. 

A key substantive difference between a conformal structure and a
projective structure a projective structure carries a well defined
notion of a round disk. Let $\Sigma$ be projective structure.  Then a {\em
  round disk} on $\Sigma$ is a projective map from a round disk in
$\chat$ to $\Sigma$. If $M$ is a hyperbolic 3-manifold then a {\em
  half-space} in $M$ is a local isometry from a half space in
$\hthree$ to $M$. Note that the projective boundary of a half-space in
$\hthree$ is a round disk so every half space in a hyperbolic 3-manifold
extends to a round disk on the projective
boundary. 

The following result is our generic inflexibility theorem for
Schwarzian derivatives. It should be compared to
Theorems~\ref{maininflex} and~\ref{mainlengthinflex}. In a future paper we will apply this result to hyperbolic cone-manifolds.
\begin{theorem}\label{mainschwarzinflex}
Let $g_t$, $t \in [a,b]$, be a one-parameter family
  of hyperbolic metrics on the interior of a 3-manifold $M$ with
  boundary.  Let $\eta_t$ be the time $t$ derivative
  of the metrics $g_t$ and let $N_t$ be a family of submanifolds of
  $M$ with compact boundary such that $\eta_t$ is a harmonic strain
  field on $N_t$. Also assume that
$$\sqrt{\int_{N_t} \|\eta_t\|^2 + \|D_t \eta_t\|^2} \leq K$$
for some $K>0$.  Let $S$ be a component of $\del M$ such that each
hyperbolic metric $g_t$ extends to a fixed conformal structure $X$ on
$S$ and a family of projective structures $\Sigma_t$ on $S$. Assume
that every embedded round disk in $\Sigma_t$ bounds an embedded
half space $H$ in $N_t$ and that
$$d_{M_t}(H, M \backslash N_t) \geq d$$
for some $d>0$.  Then
$$d(\Sigma_a, \Sigma_b) \leq CKe^{-d}$$
where $C$ is a constant depending on $\|\Sigma_a\|_F$ and the
injectivity radius of the hyperbolic metric on $X$.
\end{theorem}

\bold{Proof.} Let $H$ be an embedded half space in $M_t$ bounding a
round disk in $\Sigma_t$. By Theorem~\ref{expdecay} we have
$$\int_H \|\eta_t\|^2 + \|D_t \eta_t\|^2 \leq K^2 e^{-2d}.$$
Let $\Phi_t$ be the holomorphic quadratic differential that is the
time $t$ derivative of the family of projective structures
$\Sigma_t$. Then by Theorem 5.5 in \cite{Bromberg:Schwarzian} we have
$$K^2e^{-2d} \geq 2\sqrt{\frac{2\pi}{3}} \frac{\tanh^2 (\kappa/2)}{1 + 2\|\Sigma_t\|_F} \|\Phi_t\|_\infty$$
where $\kappa$ is the injectivity radius of the hyperbolic structure
on $X$. Integrating this inequality finishes the proof of the
theorem. For details see the proof of Theorem 1.3 in
\cite{Bromberg:Schwarzian}. \qed{mainschwarzinflex}

We will now apply this theorem to quasiconformal deformations of complete hyperbolic manifolds where some components of the conformal boundary our fixed. We will be interested in measuring the change in projective structures for these fixed components of the conformal boundary. A typical example is the deformation of a quasi-Fuchsian manifold in a Bers slice. See for example Theorem \ref{bersconv}. 

Let $M$ be a complete hyperbolic 3-manifold. Then each component
$X$ of the conformal boundary of $M$ will bound a component of
$M \backslash C(M)$, the complement of the convex core.  Label this
component $\cN(X)$ which should be thought of as a standard
neighborhood of $X$ in $M$. If $X$ is a union of components
of the projective boundary then $\cN(X)$ is the corresponding
union of components of $M \backslash C(M)$. 

If $M_t$ is a one-parameter family of complete hyperbolic structures and $X$ is a component of conformal boundary that is fixed under the deformation then the notation $\cN(X)$ does not distinguish which manifold the neighborhood lies in. In this situation we will use the projective structure on $X$ to label then end. Namely if $\Sigma_t$ is the projective boundary for $X$ in the  manifold $M_t$ then $\cN(\Sigma_t)$ is the neighborhood $\cN(X)$ in $M_t$.

Theorem~\ref{qcflow} gave us one-parameter family of hyperbolic
manifold interpolating between 
%% JFB 
the domain and range of 
a quasiconformal deformation. We will need to use this result again but
we will also need to know that the corresponding strain fields are
$L^2$ in a neighborhood of those ends of the boundary where the
deformation is conformal. For convenience we restate
Theorem~\ref{qcflow} as part of the theorem below.
\begin{theorem}\label{conformalqcflow}
  Let $\Psi: M_0 \to M_1$ be a $K$-quasiconformal deformation of the
  complete orientable hyperbolic 3-manifold $M_0$.  Then there exists
  a one-parameter family, $M_t = (M, g_t)$, $t \in [0,1]$, of
  hyperbolic metrics $g_t$ with time $t$ derivative $\eta_t$ such that
  the following holds:
\begin{enumerate}
\item The $\eta_t$ are harmonic strain fields and $\|\eta_t\|_\infty,
  \|D_t \eta_t\|_\infty \leq 3k$ where $k = \frac{1}{2} \log K$;

\item Let $\Phi_t : M_0 \to M_t$ be the identity map on $M$. Then
  $\Phi_t$ is $K^{\frac{3}{2}}$-bi-Lipschitz and $\Phi_1$ is homotopic
  to $\Psi$;

\item Let $X$ be a union of components of the conformal boundary $M_0$
  such that $\Psi$ extends to a conformal map on $X$. Then $\Phi_t$
  extends to a conformal map on $X$ for all $t$ and
$$\int_{\cN(\Phi_t(X))} \|\eta_t\|^2 + \|D_t \eta_t\|^2 < \infty.$$
\end{enumerate}
\end{theorem}

\bold{Proof.} We only need to prove (3) as (1) are (2) are exactly the
same as Theorem~\ref{qcflow}. The fact that $\Phi_t$ is conformal on
$X$ follows directly from the construction in
\cite{Reimann:visual}. To establish the $L^2$-bounds, we lift $\eta_t$
to a harmonic strain field $\tilde{\eta}_t$ on the universal cover
$\hthree$. Then $\tilde{\eta}_t$ is the visual extension of a Beltrami
differential $\mu_t$ on $\chat$. By construction, $\mu_t$ will be zero
on $\tilde{\Phi}_t(\Omega_X)$ where $\Omega_X$ is the component of the
domain of discontinuity that descends to $X$.

Let $p$ be a point in $\cN(\Phi_t(X))$. There is a unique point $q$ in
$\del C(M_t)$ that is nearest to $p$. Let $\sigma$ be the shortest
geodesic between $p$ and $q$, let $\tilde{\sigma}$ be a lift of this
geodesic to $\hthree$ and let $\tilde{p}$ and $\tilde{q}$ be the
endpoints of this geodesic which lie in the pre-images of $p$ and $q$,
respectively. Let $P$ be the hyperbolic plane in $\hthree$ that
contains $\tilde{q}$ and is perpendicular to $\tilde{\sigma}$. The
boundary of $P$ is a circle in $\chat$ that bounds a disk $D$
contained in $\tilde{\Phi}_t(\Omega_X)$. An easy calculation shows
that the in the visual measure based at $\tilde{p}$, the ratio of the
area of $D$ to the area of the the entire sphere is $\tanh
d(p,q)$. This implies that
$$\|\eta_t(p)\| = \|\tilde{\eta}_t(p)\| \leq C(1 - \tanh d(p,q)) \sim 2Ce^{-2d(p,q)}$$
where $C$ is a constant that only depends on $\|\mu_t\|_\infty$. The
area of the surface obtained by taking the locus of points in
$\cN(\Phi_t(X))$ a distance $d$ from $\del C(M_t)$ grows like
$e^{2d}$. Together these two estimates imply that the integral of
$\|\eta_t\|^2$ over $\cN(\Phi_t(X))$ is finite.

To estimate the norm of $\|D_t \eta_t\|$ we note that the lift of this
strain field is obtained by averaging $\imath \mu_t$ so the same
argument shows that it has finite $L^2$-norm on
$\cN(\Phi_t(X))$. \qed{conformalqcflow}

We can now prove the quasiconformal deformation version of our
inflexibility theorem for Schwarzian derivatives.
\begin{theorem}\label{qcschwarzinflex}
  Let $\Psi: M_0 \to M_1$ be a $K$-quasiconformal deformation of
  complete, hyperbolic 3-manifolds. Assume that the conformal boundary
  of $M_0$ is the disjoint union of two collections of components $X$
  and $Y$ and that $\Psi$ extends to a conformal map on $X$. Let
  $\Sigma_0$ be the projective structure on $X$ and $\Sigma_1$ the
  projective structure on $\Psi(X)$.  Let $d$ be the minimal distance
  between $\cN(X)$ and $\cN(Y)$ in $M_0$. Then
$$d(\Sigma_0, \Sigma_1) \leq C_0e^{-C_1 d}$$
where $C_0$ and $C_1$ depend only on $K$,  the area of the hyperbolic structure on $Y$, $\|\Sigma_0\|_F$
and the injectivity radius of the hyperbolic structure on $X$.
\end{theorem}

\bold{Proof.} We want to apply Theorem~\ref{mainschwarzinflex}. Let
$M_t$ be the one-parameter family of hyperbolic 3-manifolds given by
Theorem~\ref{conformalqcflow}. Then the submanifolds $N_t$ will be the
union of the convex cores $C(M_t)$ and the neighborhoods
$\cN(\Sigma_t)$. By Lemma~\ref{finiteconvex} the $L^2$-norm of
$\eta_t$ and $D_t \eta_t$ is finite on $C(M_t)$, and by (3) of
Theorem~\ref{conformalqcflow} these $L^2$-norms are finite on
$\cN(\Sigma_t)$. Therefore the $L^2$-norms are finite on the union
$N_t$. Just as in Lemma~\ref{finiteconvex} the boundary of $\del N_t$ will not be piecewise smooth. This can be dealt with exactly as in the proof of Lemma~\ref{finiteconvex} and we can
apply Lemma~\ref{boundaryterm} to see that
$$\int_{N_t} \|\eta_t\|^2 + \|D_t \eta_t\|^2 \leq \area(Y)9k^2$$
where $\area(Y)$ is the area of the hyperbolic structure on $Y$.

The maps $\Phi_t:M_0 \to M_t$ are $K^\frac{3}{2}$-bi-Lipschitz and
such a map between hyperbolic manifolds will take a convex set to a
$K_0$-quasi-convex set where $K_0$ depends on $K$. Applying this fact
to $\Phi^{-1}_t$ we see that the Hausdorff distance between
$\Phi_t(C(M_0))$ and $C(M_0)$ is bounded by a constant $K_1$ which
again only depends on $K$. In particular the distance between
$\cN(\Phi_t(X))$ and $\cN(\Phi_t(Y))$ is bounded below by
$d/K^{\frac{3}{2}} - K_1$.

Finally we see that if $D$ is round disk in $\chat$ bounding a half
space in $H$ then $D$ descends to an embedded disk in projective
boundary of $M_t$ if every deck transformation for $M_t$ takes $D$ off
itself. But if this is the case the same will hold for $H$ so $H$ will
descend to an embedded half space in $M_t$.

We are now in position to apply Theorem~\ref{mainschwarzinflex} to see that
$$d(\Sigma_0, \Sigma_1) \leq C_0 e^{-C_1 d}$$
where $C_0 = C\area(Y)9k^2e^{-K_1}$ with $C$ the constant from
Theorem~\ref{mainschwarzinflex} and $C_1 =
1/K^\frac{3}{2}$. \qed{qcschwarzinflex}

\bold{Remark.} If the components of $X$ are incompressible then
Nehari's Theorem \cite{Nehari:schwarzian} implies that $\|\Sigma_t\|_F
\leq 3/2$. In particular, the constants in the previous theorem will not
depend on $\|\Sigma_0\|_F$ in this case.

\medskip

\bold{Remark.} As with our previous inflexibility theorems for
quasiconformal deformations, Theorem~\ref{qcschwarzinflex} also holds
for certain hyperbolic 3-manifolds with rank-one cusps. For example if
$(M_t \cup \cN(\Phi_t(Y))) \backslash (M_t)^\delta_c$ is a compact
manifold then the proof of Theorem~\ref{qcschwarzinflex} goes through
after making the exact same modifications that were described in the
remark after the proof of Theorem~\ref{qclengthinflex}.  Manifolds
lying on the boundary of a Bers slice of a closed surface give one
important case where this condition holds.

\newcommand{\cat}{{\rm CAT}}

\section{Curves on surfaces and limits of surface groups
}

The application of inflexibility to uniformization of $3$-manifolds
fibering over the circle requires us to develop some preliminary
notions from algebraic and geometric convergence of Kleinian groups.
We emphasize that the techniques we develop treat only the case when
$S$ is closed, though many results hold more generally.  We will
assume $S$ is closed in the sequel.

\bold{Hyperbolic surfaces.} We begin by reviewing some standard facts
about hyperbolic surfaces.  A proof of the following Lemma of Bers can
be found in \cite{Buser:book:spectra}.
\begin{lemma}\label{bersconstant}
  Given a closed surface $S$ of genus $g$ there exist positive $L_g$
  and $L'_g$ such that for any hyperbolic structure $X$ on $S$ the
  following holds.
\begin{enumerate}
\item For all points $p$ in $X$ there is an essential simple closed
  curve of length at most $L_g$ that contains $p$.
\item Any simple closed curve on $X$ of length at most $L_g$ can be
  extended to pants decomposition of total length at most $L'_g$.
\end{enumerate}
\end{lemma}

We will employ the thick-thin decomposition for hyperbolic surfaces as
well as hyperbolic 3-manifolds from Lemma~\ref{Margulis}.
For surfaces, the thick-part satisfies a bounded diameter condition as
an application of Gauss-Bonnet.
\begin{lemma}\label{thickdiam}
  Each component of $X^{\ge \epsilon}$ has diameter bounded by a constant
  $D$ depending only on $\epsilon$ and $S$.
\end{lemma}
A surface $X$ is {\em $\epsilon$-thick} if $X^{\ge \epsilon} = X$.

\bold{The complex of curves.}  Given a closed surface $S$ of negative
euler characteristic, let $\calS$ denote the collection of isotopy
classes of simple closed curves on $S$.  The {\em complex of
  curves} $\calC(S)$, is a simplicial complex of dimension $3g -2$
whose vertices correspond to elements of $\calS$, and whose
$k$-simplices span collections of $k+1$ vertices whose corresponding
elements of $\calS$ can be realized disjointly on $S$.  Giving each
simplex the standard metric, we obtain a distance function
$$d_\calC \colon \calS \times \calS \to \natls.$$

A standard projection map from $\Teich(S)$ to $\cC(S)$ is readily
defined by applying the following Lemma, which is a simple application
of the Collar Lemma \cite[Thm. 4.4.6]{Buser:book:spectra} and
\cite[Lem. 2.1]{Masur:Minsky:CCI}.
\begin{lemma}
\label{intersectionbound}
  Given $L>0$ there exists $C>0$ such that if $\alpha$ and $\beta$ are
  simple closed curves on $X$ of length at most $L$ then we have
  $d_\cC(\alpha, \beta) \leq C$.
\end{lemma}
The coarse projection map
$$
\pi_\calC \colon \Teich(S) \to P(\calC^0(S))
$$
of $\Teich(S)$ to the set $P(\calC^0(S))$ of subsets of vertices of
$\calC(S)$, assigns to each $X \in \Teich(S)$ the collection of
vertices of $\calC(S)$ whose corresponding curves can be realized on
$X$ with length less than $L_g$.  By Lemma \ref{bersconstant}, the
image $\pi_\calC(X)$ is non-empty and by Lemma \ref{intersectionbound}
it has uniformly bounded diameter, so we have a coarse notion of
separation between bounded length curves on $X$ and $Y$ obtained by
taking
$$d_\calC(X,Y)  = \diam_{\calC(S)}(\pi_\calC (X), \pi_\calC(Y)).$$

\bold{Thurston's compactification.}  The elements of $\calS$ naturally
determine points in Thurston's compactification for $\Teich(S)$, the
{\em projective measured lamination space} $\pml(S)$.  Thurston showed
Teichm\"uller space can be compactified by the $(6g-7)$-sphere $\pml(S)$
to obtain a closed ball.  The action of the mapping class group
$\Mod(S)$ on $\Teich(S)$ extends to the compactification by
homeomorphisms.  Each simple closed curve $\alpha$ determines a point
in $\pml(S)$.  For further details on Thurston's construction, we
point the reader to \cite{FLP}, \cite{Imayoshi:Taniguchi:book}, or
\cite{Bonahon:currents}.

\bold{Pseudo Anosov-mapping classes.}  Those elements $\psi \in
\Mod(S)$ with positive translation distance realized at a point on the
interior of $\Teich(S)$ are known as {\em pseudo-Anosov} mapping
classes.  Their action on $\calC(S)$ is characterized by a freeness
condition: for each $\gamma \in \calC^0(S)$, we have $\gamma \not=
\psi^n(\gamma)$ for any non-zero $n$.  Thurston showed these elements
have {\em north-south dynamics} on the compactified Teichm\"uller
space: there is a unique stable lamination $[\mu^+]$ and unstable
lamination $[\mu^-]$ in $\pml(S)$ fixed by the action of $\psi$, and
for each neighborhood $U$ of $[\mu^+]$ and each $[\gamma] \in \pml(S)$
with $[\gamma] \not= [\mu^-]$, there is an $n_0$ for which
$\psi^n([\gamma])$ lies in $U$ for all $n> n_0$, and similarly for $[\mu^-]$.

\bold{Surface groups.}  We discuss two related notions of convergence
for hyperbolic 3-manifolds with the homotopy type of a surface
$S$.  A sequence $\{\rho_i\}$ of discrete, faithful representations
$$\rho_i: \pi_1(S) \to \PSL_2(\cx)$$
converges to a limit $\rho_\infty$ if $\rho_i(\gamma) \to
\rho_\infty(\gamma)$ in $\PSL_2(\cx)$ for every $\gamma \in \pi_1(S)$.
The quotient topology determined by passing to conjugacy classes is
the {\em algebraic topology}, and the set of all conjugacy classes of
discrete, faithful representations of $\pi_1(S)$ to $\PSL_2(\cx)$ with
this topology is denoted $AH(S)$.

On the level of quotient hyperbolic 3-manifolds one obtains a similar
formulation of convergence via the notion of a {\em marking} of a
hyperbolic 3-manfiold by a homotopy equivalence with $S$.  Precisely,
for each $i$ let $M_i$ be a complete hyperbolic 3-manifold and $$f_i
\colon S \to M_i$$ a homotopy equivalence. Then the marked manifolds
$\{(f_i, M_i)\}$ converge to the marked manifold $(f_\infty,
M_\infty)$ if there are lifts $\tilde{f}_i : \tilde{S} \to \tilde{M}_i
=\hthree$ such that $\tilde{f}_i$ converges to $\tilde{f}_\infty$
uniformly on compact sets.  Giving such pairs the equivalence
relation $$(f,M) \sim (g,N)$$ if there is an isometry $\phi \colon M
\to N$ so that $\phi \compos f \simeq g$, the quotient topology yields
the {\em algebraic topology} on the set $\{ [(f,M)] \}$ of equivalence
classes of marked hyperbolic 3-manifolds homotopy equivalent to $S$.
The topology is equivalent to that given above for representations via
the natural bijective holonomy relation between conjugacy classes of
discrete faithful $\PSL_2(\cx)$ representations $\rho$ of $\pi_1(S)$
and equivalence classes $[(f,M)]$.  We will also use $AH(S)$ to refer
to the collection of equivalence classes of such marked hyperbolic
$3$-manifolds with the algebraic topology.  When the meaning is clear
from context, we will also refer to a hyperbolic $3$-manifold $M$ in
$AH(S)$ assuming an implicit marking by a homotopy equivalence $f
\colon S \to M$.

As in the setting of $\Teich(S)$, the mapping class group $\Mod(S)$
acts on $AH(S)$
via {\em remarking}
$$\varphi(f,M) \mapsto (f \circ \varphi^{-1},M).$$  As a result,
we have the diagonal action $\varphi(Q(X,Y)) =
Q(\varphi(X),\varphi(Y))$ of the mapping class $\varphi \in \Mod(S)$
on quasi-Fuchsian space.

\bold{Geometric convergence.}  Let $(M_n, p_n)$ be a sequence of
hyperbolic 3-manifolds with basepoint. We say that $(M_n, p_n)$
converges {\em geometrically} to a based hyperbolic 3-manifold
$(M_\infty, p_\infty)$ if for every compact subset $K$ of
$M_\infty$ containing $p_\infty$ and every $L>1$ there exist
$L$-bilipschitz embeddings
$$
\phi_n \colon (K, p_\infty) \to
 (M_n, p_n)
$$ 
for $n$ sufficiently large. The maps $\phi_n$ are the
{\em approximating maps}. We note that this form of geometric convergence 
is often called {\em  bi-Lipschitz convergence}.

The following lemma relates geometric convergence to algebraic convergence.
\begin{lemma}\label{geomtoalg}
  Let $(M_n,p_n)$ converge to $(M_G, p_G)$ geometrically. Let $f \colon S
  \to M_G$ be a map whose image is contained in an open set $\cK$
  whose closure is compact and assume $p_G \in \cK$. Let $\phi_n: (\cK,
  p_G) \to (M_n, p_n)$ be approximating maps with bi-Lipschitz constant
  limiting to $1$, and assume that $\phi_n \circ f \colon S \to M_n$ are
  homotopy equivalences. Then $(\phi_n \circ f, M_n)$ converges to
  $(f_\infty, M_\infty)$ where $M_\infty$ is the cover of $M_G$
  induced by the subgroup $f_* (\pi_1(S))$ and $f_\infty$ is the lift
  of $f$.
\end{lemma}

 \bold{Proof.} We lift the $\phi_n$ to maps $\tilde{\phi}_n :
(\tilde{\cK}, \tilde{p}_G) \to (\hthree, \tilde{p}_n)$. Note that
$\tilde{\cK}$ is a subset of $\hthree$ and we can assume that
$\tilde{p}_G = \tilde{p}_n$ and that the derivative $D \tilde{\phi}_n$
converges to the identity on the tangent space at $\tilde{p}_G$. By
Arzela-Ascoli this sequence will be pre-compact in the compact-open
topology and since the bi-Lipschitz constant limits to $1$, every
limit will be an isometry with derivative the identity on the tangent
space at $\tilde{p}_G$. Therefore $\tilde{\phi}_n$ converges to the
identity map and the lemma follows. \qed{geomtoalg}

We would like to compare an algebraic convergence to geometric
convergence. We say that an algebraically convergent sequence
$[(f_n,M_n)] \to [(f_\infty, M_\infty)]$ converges {\em
  strongly} if the following holds. Let $(f_n,M_n)$ be
representatives such that $(f_n,M_n)$ converges to $(f_\infty,
M_\infty)$ and let $p_n =f_n(p)$ where $p$ is a point in $S$. Let
$(M_G, p_G)$ be the geometric limit of $(M_n, p_n)$. Then $[(f_n,
M_n)]$ converges to $[(f_\infty, M_\infty)]$ strongly if $(M_G, p_G) =
(M_\infty, p_\infty)$.

Note that if $(f_n, M_n)$ converges to $(M_\infty, p_\infty)$ and the
convergence is strong then the approximating maps $\phi_n$ can be chosen
such that if $\cK$ is a compact set with $f_\infty(S) \subset \cK$
then $f_n$ is homotopic to $\phi_n \circ f_\infty$.

We will use the following fundamental result of Thurston and an
improvement due to R. Evans.
\begin{theorem}[Thurston, Evans]
\label{typestrong}
Let $[\rho_n] \rightarrow [\rho]$ be a convergent sequence in $AH(S)$ and
assume that for all $\alpha \in \pi_1(S)$, if $\rho(\alpha)$ is parabolic
then $\rho_n(\alpha)$ is parabolic for all $n$. Then the
convergence is strong.
\end{theorem}

\bold{Remark.} The case when $\rho_n$ is assumed quasi-Fuchsian was
established by Thurston (see \cite{Thurston:book:GTTM}), and
generalized by Evans (\cite{Evans:persists}) to setting of general
manifolds in $AH(S)$. We will use exclusively the case when $\rho$ has
no parabolic elements in its image; the proof in this setting is
considerably easier.

\subsection{Lipschitz maps}
Let $g\colon X \to M$ be a 1-Lipschitz homotopy equivalence of a
hyperbolic surface $X$ into a hyperbolic 3-manifold $M$.  If $\alpha$
is a homotopy class of simple closed curve on $X$ then the length of
the geodesic representative of $\alpha$ on $X$ bounds from above the length of
its geodesic representative in $M$.  As a result, geometric features
of hyperbolic surfaces can be used to control the geometry of
3-manifolds (cf. \cite{Thurston:hype2}, \cite{Minsky:torus},
\cite{Minsky:CKGI}, \cite{Brock:Canary:Minsky:elc}).

Two standard constructions of such maps
 are Thurston's {\em pleated surfaces} and the related {\em
  simplicial hyperbolic surfaces}, also introduced in
\cite{Thurston:book:GTTM} and used extensively by Bonahon
\cite{Bonahon:tame} and Canary \cite{Canary:ends}.
Though we will employ both constructions, we need only their
consequences rather than the constructions themselves.
\begin{theorem}[Canary]
\label{interpolation}
Let $S$ be a closed surface and let $M \in AH(S)$.  Let $x$ and $y$ be
points in the convex core of $M$. We then have a homotopy $g_t: X_t
\to M$ with the following properties.
\begin{enumerate}
\item The family $X_t$ is a continuously varying family of
  hyperbolic metrics on $S$.
\item The maps $g_t$ are 1-Lipschitz.
\item The point $x$ lies in $g_0(X_0)$ and $y$ lies in $g_1(X_1)$.
\end{enumerate}
In particular, for any point $x$ in the convex core of $M$, there is a
1-Lipschitz map of a hyperbolic surface into $M$ whose image contains
$x$.
\end{theorem}

The previous result can be proven using simplicial hyperbolic
surfaces. For the following, one needs pleated surfaces directly.  We
use this result only in Corollary~\ref{thinbounds}.

\begin{prop}
\label{realize}
Let $\alpha^*$ be a closed geodesic in $M \in AH(S)$ that is homotopic
to a simple closed curve $\alpha$ on $S$. Then there is a 1-Lipschitz map $X
\to M$ of a hyperbolic surface $X$ that restricts to an isometry from
the geodesic representative of $\alpha$ on $X$ to $\alpha^*$.
\end{prop}
In this case we say that $X$ {\em realizes} $\alpha$.

The following lemma recapitulates a standard fact for pleated surfaces
(see \cite{Thurston:hype2}) in the setting of Lipschitz homotopy
equivalences of hyperbolic surfaces and $3$ manifolds.  It will be
useful to know this for arbitrary Lipschitz constants.
\begin{lemma}
\label{catthin}
Given $\epsilon>0$ and $B\ge 1$ there exists $\epsilon'>0$ such that
if $f: X \to M$ is a $B$-Lipschitz homotopy equivalence of a
hyperbolic surface into a hyperbolic 3-manifold $M$ and $p$ is a point
with $f(p) \in M^{< \epsilon'}$ then we have $p \in X^{<
  \epsilon}$.
\end{lemma}

 \bold{Proof.} 
By Lemma \ref{thickdiam}, the diameter of each component of $X^{\ge
  \epsilon}$ is bounded by a constant $D$ that only depends on $S$ and
$\epsilon$.  Therefore the $f$-image of each component of $X^{\ge
  \epsilon}$ has diameter less than $BD$. By a theorem of Brooks and
Matelski (see \cite{Brooks:Matelski:collars}), we may choose
$\epsilon' < \varepsilon_3$ small enough such that the distance
between the boundaries of the $\varepsilon_3$-thin and the $\epsilon'$-thin
part is at least $BD$. Every component of $X^{\ge \epsilon}$ has
non-abelian fundamental group while every component of $M^{<
  \varepsilon_3}$ has abelian fundamental group.  Since $f$ is
$\pi_1$-injective, the $f$-image of each component of $X^{\ge \epsilon}$
must intersect $M^{\ge \varepsilon_3}$ and is therefore disjoint from
$M^{< \epsilon'}$.
\qed{catthin}

Mumford's compactness theorem (see \cite{Mumford:compact}) guarantees
that any sequence of $\epsilon$-thick surfaces in $\Teich(S)$ can be
re-marked to converge in $\Teich(S)$ up to subsequence; the following
shows the same is true for $M_n \in AH(S)$ with uniformly Lipschitz
markings by thick surfaces.
\begin{prop}
 Let $\epsilon >0$ be given.
\begin{enumerate} 
\item  For each sequence $\{X_n\}$ of $\epsilon$-thick surfaces there are
  markings $f_n \colon S \to X_n$ such that $(f_n , X_n)$ converges in
  $\Teich(S)$.
\item Let $(f_n,X_n)$ be a convergent sequence in $\Teich(S)$ and $g_n
  \colon X_n \to M_n$ $B$-Lipschitz homotopy equivalences to
  hyperbolic $3$-manifolds $M_n$. Then $\{(g_n \circ f_n, M_n)\}$ has
  a convergent subsequence in $AH(S)$.
\end{enumerate}
\label{algebraiclimit}
\end{prop}

 \bold{Proof.} Statement (1) is a restatement of Mumford's
compactness theorem for the Moduli space $\calM(S)$
\cite{Mumford:compact}.  To see statement (2), note that since the sequence
$( f_n, X_n)$ converges, we can place a hyperbolic metric on $S$ such that
the marking maps $f_n$ are $B'$-Lipschitz for some $B'>1$. Then the
maps $h_n = g_n \circ f_n$ are $BB'$-Lipschitz.

Pick a point $p \in S$ and let $\tilde{p} \in \tilde{S} = \htwo$ be a
point in the pre-imiage of $p$. Identifying each $\tilde{M}_n$ with
$\hthree$, we choose lifts of $h_n$ such that $\tilde{h}_n(p) = 0 \in
\hthree$.  Since the maps $h_n$ are $BB'$-Lipschitz, it follows that
for any $q \in \tilde{S}$, the set $\{\tilde{h}_n(q)\}$ has compact
closure in $\hthree$.  By the Arzela-Ascoli theorem, there exists a
subsequence such that $\tilde{h}_n$ converges uniformly on compact
sets to a map $\tilde{h}_\infty: \tilde{S} \to \hthree.$ The action of
$\pi_1(S)$ on $\tilde{S}$ commutes with the action of a representation
of $\pi_1(S)$ in $\PSL_2(\cx)$ so that $\tilde{h}_\infty$ descends to
a pair $(h_\infty, M_\infty)$ where $M_\infty$ is the quotient
$3$-manifold.  \qed{algebraiclimit}

\subsection{Margulis estimates}
Let $M$ be a hyperbolic manifold in $AH(S)$ and $K$ a subset of $M$.
The Margulis lemma provides bounds for the number of homotopy classes
of essential primitive loops of length less than $L$ that intersect
$K$ such that each loop is a homotopic to a simple closed curve on
$S$.

\begin{lemma}
\label{margulisbound}
Given $L>0$ and $D>0$ there is a $N>0$ such that the following
holds. Let $M \in AH(S)$ and let $K \subset M$ be a subset of
diameter at most $D$.  Then the number of distinct essential homotopy
classes of loops of length at most $L$ intersecting $K$ is bounded
above by $N$.
\end{lemma}

 \bold{Proof.}
By \cite{Brooks:Matelski:collars} there exists $\epsilon>0$ such
that the distance between $\del M^{< \epsilon}$ and $\del M^{\ge
  \varepsilon_3}$ is at least $D+L$ where $\epsilon$ is less than the
3-dimensional Margulis constant $\varepsilon_3$.

The proof then breaks into two cases. First assume that $K$
intersects $M^{< \epsilon}$. Then every loop of length at most $L$
that intersects $K$ will be contained in a component of $M^{<
  \varepsilon_3}$.  Since $M$ lies in $AH(S)$, $M$ has no rank-two
cusps and every component of $M^{< \varepsilon_3}$ contains one
essential primitive loop.
% so we have $\#(K, L) \leq 1$.

Now we assume that there is a point $x$ in $K \cap M^{\ge \epsilon}$.
Any loop of length at most $L$ that intersects $K$ will be homotopic
to a loop of length at most $L +2D$ that intersects $x$. The number of
distinct homotopy classes of loops of length at most $L + 2D$ that
intersect $x$ is bounded by the quotient $$V = \frac{\vol
  (B_{\half^3}(0,L+2D+\epsilon))}{\vol (B_{\half^3}(0,\epsilon))}$$ of
the volumes of balls of radius $L+2D + \epsilon$ and $\epsilon$ about
the origin in hyperbolic space $\half^3$, so taking $N = V +1$ proves
the Lemma.  \qed{margulisbound}

\newcommand{\short}{{\bf short}}
\subsection{Geometric limit arguments}
A subset $K \subset M$ is {\em $\epsilon$-thick} if $K \subset M^{\ge
  \epsilon}$.

\begin{prop}
\label{thickstrong}
Let $(M_n, \omega_n)$ be hyperbolic 3-manifolds homotopy equivalent to
$S$ that converge geometrically to $(M_\infty, \omega_\infty)$. Assume
there exist $\epsilon > 0$ and $R_n \to \infty$ such that the
$R_n$-neighborhood of $\omega_n$ in $M_n$ is $\epsilon$-thick.  Then
$M_\infty$ is homotopy equivalent to $S$, and there are
homtopy equivalences $f_n \colon S \to M_n$ and $f_\infty
\colon S \to M_\infty$ so that $(f_n, M_n)$ converges strongly to
$(f_\infty, M_\infty)$.
  \marginpar{\tiny We need to assume that
  $M$ has no parabolics if we only want to use Theorem
  \ref{typestrong}. This is good enough for the application to double
  limits. I guess we could also uses Richard's theorem. K}
\end{prop}

 \bold{Proof.} By Theorem \ref{interpolation} there is a
hyperbolic surface $X_n$ and a 1-Lipschitz map $g_n: X_n \to M_n$
whose image contains $\omega_n$. Let $q_n$ be a point in $X_n$ with
$g_n(q_n) = \omega_n$. Since $g_n$ is 1-Lipschitz, an
$R_n$-neighborhood of $q_n$ in $X_n$ will be $\epsilon$-thick as
well. There is a constant $K$ depending only on $\epsilon$ and the
genus of $S$ such that if a hyperbolic structure $X$ on $S$ has
$\epsilon$-thick neighborhood of radius at least $K$ then $X$ itself
is $\epsilon$-thick. In particular for large $n$ the surfaces $X_n$
are themselves $\epsilon$-thick.

We now apply Proposition \ref{algebraiclimit} to find homeomorphisms
$f_n \colon S \to X_n$ such that $\{(f_n, X_n)\}$ converges in
Teichm\"uller space and $\{(g_n \circ f_n, M_n)\}$ converges in
$AH(S)$.  To show the sequence converges strongly it suffices to
verify it is type-preserving by an application of
Theorem~\ref{typestrong}.  After an isotopy, we can assume there is a
fixed point $x \in S$ such that $f_n(x) = q_n$. Let $\alpha$ be a
non-trivial loop in $S$ based at $x$. Since $\{( f_n,X_n)\}$ converges
we can homotope the $f_n$ rel $x$ so that the loops $f_n(\alpha)$ have
length bounded above by a constant only depending on the homotopy
class of $\alpha$ rel $x$.  Since each $g_n$ is 1-Lipschitz, the
lengths of the loops $g_n \circ f_n(\alpha)$ are also uniformly
bounded.  If the sequence is not type-preserving there will be some
$\alpha$ such that the length of the geodesic representative of $g_n
\circ f_n(\alpha)$ tends to zero.  In particular, for large $n$ the
curve $g_n \circ f_n(\alpha)$ will be homotopic into a component of
the $\epsilon$-thin part of $M_n$.  \marginpar{\tiny we use this fact
  multiple times -- make a lemma? -J} There is then a bound on the
distance from $g_n \circ f_n(\alpha)$ to this component of the thin
part where the bound only depends on the length of $g_n \circ
f_n(\alpha)$.  But for large $n$ the $R_n$-neighborhood $\calN_{R_n}(\omega_n)$
of $\omega_n$ has non-empty intersection with this component of the
$\epsilon$-thin part, contradicting our assumption that
$\calN_{R_n}(\omega_n)$ lies in $M_n^{\ge \epsilon}$.

If follows that the sequence $( g_n \circ f_n, M_n)$ is
type-preserving, and by Theorem~\ref{typestrong} the convergence is
strong. The proposition then follows. \qed{thickstrong}

\marginpar{\tiny  proof of last statement of the theorem? -J}

\begin{prop}
  Given positive constants $L$ and $\epsilon$, there exist $R$ and $C$
  so that the following holds. Let $M \in AH(S)$, and $\alpha$ and
  $\beta$ curves in $\calC(S)$. Let $\alpha^*$ and $\beta^*$ be loops
  based at $\omega$ in the convex core of $M$ in the homotopy class of
  $\alpha$ and $\beta$, respectively, and assume that
  $\ell_M(\alpha^*) \le L$, $\ell_M(\beta^*) \le L$ and the
  neighborhood $\calN_R(\omega)$ of radius $R$ about $\omega$ has
  injectivity radius bounded below by $\epsilon$. Then we have 
$$d_\calC(\alpha,\beta) \le C.$$
\label{boundedgeometrylimit}
\end{prop}

 \bold{Proof.} We argue by contradiction. Assume there is a
sequence $(M_n, \omega_n)$ of hyperbolic manifolds with baseframes
such that $\calN_{R_n}(\omega_n)$ is $\epsilon$-thick, $R_n \to
\infty$, and that $\alpha_n$ and $\beta_n$ are homotopy classes in
$\calC(S)$ represented by closed loops $\alpha^*_n$ and $\beta^*_n$ in
$M_n$ based at $\omega_n$ of length at most $L$ for which
$d_\calC(\alpha_n, \beta_n) \to \infty$.

After passing to a subsequence, $(M_n, \omega_n)$ converges
geometrically to a manifold $(M_\infty, \omega_\infty)$. By
Proposition \ref{thickstrong}, $(M_\infty, \omega_\infty)$ is homotopy
equivalent to $S$, and the approximating maps are homotopy
equivalences for large $n$.  Choosing a compact core $K$ of $M_\infty$
that contains a diameter $4L$ neighborhood of $\omega_\infty$, there
are 2-bi-Lipschitz approximating maps $\phi_n \colon K \to M_n$ for large
$n$ such that $\phi_n$ are homotopy equivalences. 

The image of $K$ under $\phi_n$ will contain $\alpha^*_n$ and
$\beta^*_n$ so $\phi^{-1}_n(\alpha^*_n)$ and $\phi^{-1}_n(\beta^*_n)$
are loops in $M_\infty$ of length at most $2L$. Since $K$ is compact,
there are only finitely many free homotopy classes of loops in $K$ of
length at most $2L$.  This finite set of loops has finite diameter in
$\calC(S)$.  Since $\phi_n$ is a homotopy equivalence, we conclude
$d_{\calC}(\alpha_n, \beta_n)$ is uniformly bounded, contrary to our
assumption.  \qed{boundedgeometrylimit}

Given $(f,M) \in AH(S)$, and $\epsilon >0$, we let
$\short_\epsilon(M)$ denote the set of isotopy classes $\alpha \in
\calC(S)$ so that $\ell_M(\alpha) < 2\epsilon$.  It is due to
Thurston, and a consequence of Lemma~\ref{catthin} and
Theorem~\ref{interpolation}, that there is an $\epsilon_{\rm s} >0$ so
that for $\epsilon < \epsilon_{\rm s}$, a closed geodesic in $M$ of
length less than $\epsilon$ lies in the homotopy class of a simple
closed curve on $S$.  Then for $\epsilon< \epsilon_{\rm s}$ and for
each $\alpha \in \short_\epsilon(M)$ there is a component
$\tube_{\epsilon}(\alpha)$ of the $\epsilon$-thin part of $M$.  
We record the following immediate consequence.
\begin{lemma}
\label{shortemptythick}
  Given $(f,M) \in AH(S)$, and positive $\epsilon < \epsilon_{\rm s}$,
  if $M$ has no cusps, and $\short_\epsilon(M) = \emptyset$, then $M$
  is $\epsilon$-thick.
\end{lemma} 

% \bold{Proof:}
% On the contrary, if $M$ were not $\epsilon$ thick, then the absence of
% parabolics implies the existence of a component
% $\tube_{\epsilon}(\gamma)$ of the $\epsilon$-thin part of $M$, a tubular
% neighborhood of a geodesic $\gamma^*$ of length less than $2\epsilon$.  Since
% $\epsilon< \epsilon_{\rm s}$, we have a simple curve $\gamma \subset S$ so
% that $f (\gamma)$ is homotopic to $\gamma^*$, so $\short_\epsilon(M) \not=
% \emptyset.$
% \qed{shortemptythick}

\begin{prop}
\label{thickballs}
Given positive $\epsilon < \epsilon_{\rm s}$, and $R>0$
there exists $L>0$ so that the following holds.  Let $f \colon X \to
M$ be a 1-Lipschitz homotopy equivalence of a hyperbolic surface $X$
into a hyperbolic $3$-manifold $M$, such that $X$ is $\epsilon$-thick.
If each $\gamma \in \short_\epsilon(M)$ satisfies $\ell_X(\gamma) > L$
then the $R$-neighborhood $\calN_R(f(X))$ is $\epsilon$-thick.
\end{prop}

% >>>PROBABLY SOME DECOUPLING OF CONSTANTS NEEDED HERE>>>

\marginpar{\tiny Todo -- fix $\epsilon$ vs. $2\epsilon$ discrepancy}
\bold{Proof.} Again we argue by contradiction and assume that we have
a sequence $\{g_n \colon X_n \to M_n\}$ of 1-Lipschitz homotopy
equivalences from $\epsilon$-thick surfaces $X_n$ with the property
that the infimum of $\ell_{X_n}(\gamma)$ over all $\gamma \in
\short_\epsilon(M_n)$ is at least $L_n \to \infty$, but the $R$-ball
about $g_n(X_n)$ is not $\epsilon$-thick for any $n$.

By Proposition \ref{algebraiclimit}, there are markings $f_n \colon S
\to X_n$ so that after passing to a subsequence $\{(f_n,X_n)\}$
converges in $\Teich(S)$ and $\{(g_n \compos f_n , M_n)\}$ converges in
$AH(S)$ to an algebraic limit $(g_\infty \compos f_\infty,M_\infty)$
with the property that $\short_\epsilon(M_\infty) = \nullset$.  Otherwise
there is a  $\gamma$  in $\short_\epsilon(M_n)$ for sufficiently large $n$, so we have
$$
\ell_{X_n}(\gamma) \rightarrow \infty.
$$ 
On the other hand, convergence of $\{(f_n,X_n)\}$ in $\Teich(S)$
implies $\ell_{X_n}(\gamma)$ converges, a contradiction.

By Theorem~\ref{typestrong}, the sequence $\{(g_n \compos f_n ,
M_n)\}$ converges strongly to $(g_\infty \compos f_\infty,M_\infty)$
with $\short_\epsilon(M_\infty) = \nullset$; in particular, by
Lemma~\ref{shortemptythick}, $M_\infty$ is $\epsilon$-thick.  By
geometric convergence, the $R$-neighborhood about $g_n(X_n)$ is
$\epsilon$-thick for $n$ sufficiently large, a contradiction completing
the proof.  \qed{thickballs}

\begin{cor} Given positive $R$, $L$ and $\epsilon < \varepsilon_3$ there are $C$ and $D$ so
  that the following holds: let $\alpha^*$ be a loop in a
  manifold $M \in AH(S)$ in the homotopy class of $\alpha \in
  \calC(S)$. Assume that length of $\alpha^*$ is at most $L$ and that
  the $R$-neighborhood of $\alpha^*$ is not $\epsilon$-thick. Then
  there is a curve $\beta \in \short_\epsilon(M)$ satisfying
$$d_\calC(\alpha,\beta) <C$$ with the property that
$d_M(\tube(\alpha),\tube(\beta)) \le D$.
\label{thinbounds}
\end{cor}

\bold{Proof.} If $\alpha^*$ is not a geodesic then it is
either uniformly close to its geodesic representative \marginpar{\tiny
  again this ruled annulus business -J} or $\alpha$ is in
$\short_\epsilon(M)$. In the latter case, we may take $\beta = \alpha$
and we are done.  Thus we can assume that $\alpha^*$ is a geodesic and
consider 1-Lipschitz hyperbolic surface $f \colon X \to M$ realizing
$\alpha^*$.

If the surface $X$ fails to be $\epsilon$-thick itself, then the
theorem follows trivially from Lemmas \ref{bersconstant} and
\ref{catthin}.  Thus we may assume that $X$ is $\epsilon$-thick.

Applying Proposition~\ref{thickballs}, given $R$ there is an $L'$ so
that if the $R$-ball about $X$ fails to be $\epsilon$-thick there
is a curve $\beta \in \calC(S)$ so that $\ell_X(\beta) < L'$.  Since
$\alpha$ has length at most $L$ on $X$, by Lemma \ref{bersconstant}
there is a $C$ depending on $\max\{L,L'\}$ with the property that
$$d_\calC(\alpha,\beta) < C.$$  

Since $X$ is itself $\epsilon$-thick, there is a uniform bound
depending only on $\epsilon$ and the genus of $X$ for the diameter of
$X$.  Hence there is a uniform bound on the distance between the
geodesic representatives of $\alpha$ and $\beta$ on $X$. For any loop
$\gamma$ in $M$ of length at most $\max\{L, L'\}$ there is a bound,
depending only on $\max\{L, L'\}$, on $d_M(\gamma,
\tube(\gamma))$. Combining the two bounds gives the result.
\qed{thinbounds}

We can now prove the main theorem of this section, providing a linear
lower bound on the distance between two bounded length curves in a
hyperbolic manifold in terms of the distance of the curves in the
curve complex. By the Margulis lemma, a short curve will have a large
tubular neighborhood and therefore lie at large distance from the
geodesic representatives of every other
bounded length curve. In this
case we will prove a stronger statement and bound the distance between
the Margulis tubes. For this reason, we define
$\tube'_\epsilon(\gamma) = \tube_\epsilon(\gamma)$ if $\ell_M(\gamma)
< \epsilon$ and let $\tube'_\epsilon(\gamma)$ be the geodesic
representative of $\gamma$ in $M$ if $\ell_M(\gamma) \geq \epsilon$.
\begin{theorem}
\label{distanceboundscurvecomplex}
Given $L>0$ there exist $K_1$ and $K_2$ all positive so that
for  $\alpha$ and $\beta$ in $\calC^0(S)$, and $M \in AH(S)$, the
following holds:  if $\ell_M(\alpha) < L$ and $\ell_M(\beta) < L$,
then 
$$d_M(\tube'_{\varepsilon_3}(\alpha),\tube'_{\varepsilon_3}(\beta)) \ge K_1
d_\cC(\alpha,\beta) - K_2.$$
\end{theorem}

\bold{Remark.}  We point out that
Theorem~\ref{distanceboundscurvecomplex} uses in an essential way the
fact that $S$ is a closed surface.  If $S$ has boundary, the same
statement holds if we measure distance in the {\em pared manifold}
$M^0$ obtained by excising cusps associated to $\bdry S$.  All the
results of the paper would then generalize in the presence of the
appropriate generalization of the geometric inflexibility theorem
(Theorem~\ref{thickinflex}) to this pared setting.

\medskip

Before we begin the proof of Theorem \ref{distanceboundscurvecomplex}
we make a definition and prove a preliminary lemma. A {\em $D$-coarse
  path} in $\cC(S)$ is a sequence of $\alpha_i$ in $\cC^0(S)$ such
that $d_\cC(\alpha_i, \alpha_{i+1})\leq D$.

\begin{lemma}\label{coarsepath}
  Given $L>0$ there exists a $D>0$ and $R>0$ such that the following
  holds. Let $\alpha$ and $\beta$ in $\cC(S)$ and $M \in AH(S)$
  satisfy $\ell_M(\alpha) \le L$ and $\ell_M(\beta) \le L$. Let $\Gamma$ be
  a path in $M$ from $\tube'_{\varepsilon_3}(\alpha)$ to
  $\tube'_{\varepsilon_3}(\beta)$. Then there are closed curves
  $\alpha_i$ with $\ell_M(\alpha_i) \leq L_g$ and
  $d_M(\alpha_i,\Gamma)<R$ such that the curves $\alpha_i$ describe a
  { $D$-coarse path} in $\cC(S)$ from $\alpha$ to $\beta$.
\end{lemma}

 \bold{Proof.} Let $x$ be the endpoint of $\Gamma$ on
$\tube'_{\varepsilon_3}(\alpha)$ and $y$ the endpoint of $\Gamma$ on
$\tube'_{\varepsilon_3}(\beta)$. Let $g_t: Z_t \to M$, $t\in [0,1]$ be
a continuous family of 1-Lipschitz maps of hyperbolic surfaces $Z_t$
such that $x \in g_0(Z_0)$ and $y \in g_1(Z_1)$. Such an interpolation
exists by Theorem \ref{interpolation}.

There is a subinterval $[a,b] \subseteq
[0,1]$ such that $x \in g_a(Z_a)$, $y \in g_b(Z_b)$ and $g_t(Z_t) \cap
\Gamma \not= \nullset$ for all $t \in [a,b]$. Reparameterize $[a,b]$
to be the interval $[0,1]$ and replace the original homotopy with this
reparametrized homotopy.

Given a simple closed curve $\gamma$ on $S$ let $U(\gamma) \subseteq
[0,1]$ be the set of $t$ such that there is a simple closed curve
$\gamma'$ on $Z_t$, homotopic to $\gamma$, with $\gamma' \cap
g_t^{-1}(\Gamma) \neq \nullset$ and $\ell_{g_t} (\gamma') < L_g$. By
(2) of Lemma~\ref{bersconstant}, if $U(\gamma) \cap U(\gamma') \neq
\nullset$ then 
\begin{equation}
d_{\cC}(\gamma, \gamma') \leq C.
\end{equation}
\marginpar{\tiny this is an application of the collar lemma right?}

Let $z$ be a point in $g_t^{-1}(\Gamma)$. By Lemma
\ref{bersconstant}  for each $t$
there exists $\gamma \in \cC(S)$ such that $t \in U(\gamma)$.  The
open (possibly disconnected) subsets $U(\gamma)$ cover $[0,1]$ so we can find a 
collection $\alpha_0, \dots, \alpha_n$ of distinct homotopy classes of
simple closed curves in $\calC(S)$ such that the $U(\alpha_i)$
satisfy
\begin{equation}
\label{intervalscover}
U(\alpha_i) \cap U(\alpha_{i+1}) \neq \nullset
\end{equation}
with $0 \in U(\alpha_0)$ and $1 \in U(\alpha_n)$. In particular the
$\alpha_i$ are a $C$-coarse path. To finish the proof we need to show
that $\alpha$ and $\beta$ are uniformly close to $\alpha_0$ and
$\alpha_n$, respectively.

To see this we observe that if $\ell_M(\alpha)$ is sufficiently small
then Lemma~\ref{catthin} guarantees that if $0 \in U(\gamma)$ we have
$\gamma = \alpha$. On the other hand if $\alpha$ has a sufficiently
large thick neighborhood then since $\ell_M(\alpha) \le L$ Proposition
\ref{boundedgeometrylimit} implies that if $0 \in U(\gamma)$ then
$\alpha$ and $\gamma$ are uniformly close in $\cC(S)$. If neither of
these cases holds, an application of Corollary~\ref{thinbounds} allows
us to replace $\alpha$ with a curve $\alpha'$ that is sufficiently
short \marginpar{\tiny what does this mean?} so that $\alpha$ and
$\alpha'$ are uniformly close in $M$ and their corresponding vertices
are uniformly close in $\cC(S)$. We then append to $\Gamma$ a geodesic
segment of length at most $ R$ connecting $x$ to
$\tube'_{\varepsilon_3}(\alpha')$ to make a new path $\Gamma'$ and
apply the previous argument to $\Gamma'$. This process yields a coarse
path $\{\alpha_i\}$ with $\alpha_0 = \alpha'$ such that the $\alpha_i$
have representatives in $M$ of length at most $L_g$ and so that each
$\alpha_i$ intersects $\Gamma'$.

Applying the same analysis to $\beta$ we obtain the desired  coarse
path.
\qed{coarsepath}

\bold{Remark.} By the Bers inequality (see \cite[Thm. 3]{Bers:bdry},
\cite[Prop. 6.4]{McMullen:iter}) given the quasi-Fuchsian manifold
$Q(X,Y)$ we have $\ell_{Q(X,Y)}(\gamma) \leq 2\ell_X(\gamma)$.
Therefore there is a uniform $D$, depending only on $S$, such that
there is a $D$-coarse path from $X$ to $Y$ in $\cC(S)$ whose lengths
in $Q(X,Y)$ are at most $L_g$.  \medskip

\bold{Proof of \ref{distanceboundscurvecomplex}.} Let $\Gamma$ be the
shortest geodesic from $\tube'_{\varepsilon_3}(\alpha)$ to
$\tube'_{\varepsilon_3}(\beta)$ so that
$d_M(\tube'_{\varepsilon_3}(\alpha), \tube'_{\varepsilon_3}(\beta)) =
\ell(\Gamma)$.  Applying Lemma~\ref{coarsepath}, with $L$ as in the
statement, we obtain a $D> 0$, an $R> 0$, and a $D$-coarse path $\{\alpha_i\}_{i=0}^n
\subset \calC^0(S)$ with $\alpha_0 = \alpha$ and $\alpha_n =\beta$,
with each $\alpha_i$ satisfying $\ell_M(\alpha_i) < L_g$ for $0<i<n$
and with each having distance at most $R$ from $\Gamma$.
Since this path $\{\alpha_i\}$ is $D$-coarse we have $n \ge
d_\calC(\alpha,\beta)/D$. 

Since $\alpha_i$ has length at most $L_g$ for $0 < i < n$ and lies at
distance at most $R$ from $\Gamma$, each determines a curve of length
at most $2R + L_g$ that intersects $\Gamma$ in the same homotopy
class.  The number of homotopically distinct primitive closed loops of
length at most $L_g+2R$ intersecting $\Gamma$ such that each is
homotopic to an essential simple curve on $S$ is at least $n$, where
$D(n+1) \geq d_{\cC}(\alpha, \beta)$.
%%%
%%% Have a look at this again
%%%
\marginpar{\tiny Have a look at this again}

The path $\Gamma$ can be divided into $\lfloor \ell(\Gamma) \rfloor$
disjoint segments of length $1$ and one segment of length at most
$1$. Let $N$ be the constant given by Lemma \ref{margulisbound} for
the length bound $L_g + 2R$ and the diameter bound $1$. Then each of
the $\lfloor \ell(\Gamma) \rfloor + 1 = \lceil \ell(\Gamma) \rceil$
segments intersects at most $N$ distinct homotopy classes of closed
curves of length at most $L_g + 2R$ and therefore
$$N(\ell(\Gamma) + 1)\ge N\lceil \ell(\Gamma) \rceil \ge n.$$
Combining this inequality with a lower bound on $n$ we have
$$N(\ell(\Gamma) + 1) \geq \frac{d_{\cC}(\alpha, \beta)}{D} - 1$$
as desired.
\qed{distanceboundscurvecomplex}

\bold{Remark.} Note that if $\alpha$ is a closed curve in $M$ of length
at most $L$ then the distance between $\alpha$ and
$\tube'_{\varepsilon_3}(\alpha)$ is uniformly bounded by a constant
only depending on $L$. In particular Theorem
\ref{distanceboundscurvecomplex} holds if we replace
$\tube'_{\varepsilon_3}(\alpha)$ with any curve of length at most $L$
that is homotopic to $\alpha$.

\begin{cor}
\label{quasifuchsianthick}
Given a closed surface $S$, there is linear function $f$ such that the
distance between the boundary components of the convex core
$C(Q(X,Y))$ of a quasifuchsian manifold $Q(X,Y)$ in $QF(S)$ is bounded
below by $f(d_{\cC}(X,Y))$.
\end{cor}

 \bold{Proof.} Let $\alpha \in \cC(S)$ have length on $X$ of
at most the Bers constant $L_g$, and choose $\beta \in \cC(S)$
similarly for $Y$. Let $\alpha^*$ and $\beta^*$ be the geodesic
representatives of $\alpha$ and $\beta$, respectively, in $Q(X,Y)$.
By the Bers inequality $\alpha^*$ and $\beta^*$ have length at most
$2L_g$.  Every closed geodesic is contained in the convex core
$C(Q(X,Y))$ so there are hyperbolic surfaces $Z_\alpha$ and $Z_\beta$
and 1-Lipschitz homotopy equivalences, $f_\alpha \colon Z_\alpha \to
Q(X,Y)$ and $f_\beta \colon Z_\beta \to Q(X,Y)$, realizing $\alpha^*$
and $\beta^*$, respectively. These maps will have image in the convex
core and each separates the two components $\bdry C(Q(X,Y))$.

Let $\Gamma$ be the shortest geodesic connecting the two components of
the convex core boundary.  Let $x$ be a point of intersection of
$\Gamma$ with the image of $Z_\alpha$ and $y$ a point of intersection
of $\Gamma$ with the image of $Z_\beta$. There is a curve $\alpha_0
\in \cC(S)$ such that $\alpha_0$ has a representative on $Z_\alpha$ of
length at most $L_g$ and whose image $\alpha^*_0$ intersects $x$. Note
that $\alpha^*_0$ will have length at most $L_g$ constant and, by
(2) of Lemma \ref{bersconstant}, there is a constant $D$ such that
$d_{\cC(S)}(\alpha, \alpha_0) \le D$. 
Similarly, we can find a curve
$\beta_0 \in \cC(S)$ that is represented by a loop $\beta^*_0$ that
intersects $y \in Q(X,Y)$, of length at most the Bers constant and
with $d_{\cC}(\beta, \beta_0) \le D$.

We want to find a lower bound for 
$\ell(\Gamma)$. We observe that $$d_{Q(X,Y)}(\alpha^*_0, \beta^*_0) \le
\ell(\Gamma) + L_g$$ and that $$d_{\cC(S)}(X,Y) =
d_{\cC(S)}(\alpha, \beta) \le d_{\cC(S)}(\alpha_0, \beta_0) + 2D.$$ The
result then follows from Theorem
\ref{distanceboundscurvecomplex}. \qed{quasifuchsianthick}

The following Corollary controls the depth of a given curve in the
convex core.  We leave the proof, which follows the same lines as the
above, to the reader.

\begin{cor}
\label{curvedeepincore}
Given a closed surface $S$ and $L>0$ there is a linear function $f_L$
such that if $\gamma \in \cC(S)$ and the length of $\gamma$ is at most
$L$ then the distance from the geodesic representative of $\gamma$ in
$Q(X,Y)$ to the boundary of the convex core is at least
$$f_L(\min\{d_{\cC}(X, \tube'_{\varepsilon_3}(\gamma)), d_{\cC}(Y, \tube'_{\varepsilon_3}(\gamma))\}).$$
\end{cor}

\section{Pseudo-Anosov double limits}
In this section, we employ the estimates on depth in the convex core
from the previous section together with the inflexibility theorems of
previous sections to establish the convergence of pseudo-Anosov double
iteration on quasi-Fuchsian space of a closed surface.  The
hyperbolization theorem for 3-manifolds that fiber over the circle
follows as a consequence.

The following Theorem is a refined version of a Theorem of Masur and
Minsky (see \cite[Prop. 3.6]{Masur:Minsky:CCI}). The proof that
follows was communicated to us by M. Bestvina who attributes the
argument to F. Luo. (See \cite[Prop.~11]{Bestvina:Fujiwara:bounded}).
\begin{theorem}
\label{soutoAnosov}
Let $\psi \in \Mod(S)$ be pseudo-Anosov with $[\mu^+]$ and $[\mu^-]$
the attracting and repelling laminations in $\pml(S)$.  Then there is
a $K_\psi$ depending only on $\psi$ so that for any $B,D>0$ the following holds.

\begin{enumerate} 

\item There exist
neighborhoods $V$ and $U$ of $[\mu^+]$ and $[\mu^-]$ in $\pml(S)$ so
that for any $\alpha \in V$ and $\beta \in U$ we
have $$d_{\cC}(\alpha, \psi^n(\beta)) \ge K_\psi n +B.$$
\item There exists a subset $W$ in $\pml(S) \backslash (U \cup V)$
  such that any path in $\cC(S)$ from a curve in $U$ to a curve in $V$
  contains a subpath of length at least $D$.

\item  For each $\alpha \in U$,
$\beta \in V$ and $\gamma \in W$ we have 
$$d_{\cC}(\gamma, \psi^{-n} (\alpha)) \ge K_\psi n + B
\ \ \ \text{and} \ \ \
d_{\cC}(\gamma, \psi^n (\beta)) \ge K_\psi n +B.$$

\item Furthermore, given any curve $\beta$, the sets $U$, $V$, and $W$
  may be taken so that any one of them contains $\beta$.
  \marginpar{\tiny We can choose which one of the three sets $\beta$
    lies in. How should we phrase this? k -- is this okay??? -J}

\end{enumerate}
\end{theorem}

 \bold{Proof.} Let $V'$ and $U'$ be neighborhoods of $[\mu^+]$
and $[\mu^-]$ in $\pml(S)$ such that for simple closed curves $\alpha
\in V'$ and $\beta' \in U'$ the intersection of $\alpha$ and $\beta$
is non-empty.  Then $V'$ and $U'$ are necessarily disjoint, and we let
$$W' = \pml(S) \setminus (V' \cup U')$$ be their complement in
$\pml(S)$. Then the north-south dynamics of $\psi$\marginpar{\tiny
  reference} guarantees that for any compact subset $\calK$ of $\pml(S)$
that does not contain $[\mu^+]$ we have $\psi^i(U')$ contains $\calK$ for
some positive $i$. In particular, since $W' \cup V'$ is compact there
is an $N$ so that we have $\psi^{N}(W' \cup V') \subset V'$.

We claim
that for any curves $\alpha \in \psi^{mN}(V')$  and $\beta \in U'$ we have
\begin{equation}
d_{\cC}(\alpha, \beta) \geq m+1.
\label{m+1}
\end{equation}

We first note that if $\alpha \in \psi^i(V')$ and $\beta \in
\psi^i(U')$ then any geodesic in the curve complex connecting them
will contain a curve in $\psi^i(W')$, the complement of the two sets,
and the distance between $\alpha$ and $\beta$ will be greater than
$2$. We also note that $\psi^{(m+1)N}(W') \subset \psi^{mN}(V')$.

We prove the inequality by induction. By the observation above, the
inequality~(\ref{m+1}) is true when $m=1$.  
To complete the induction, note that we have $\psi^i(U') \subset
\psi^{i+N}(U')$. Therefore $U'$ lies in $\psi^{(m+1)N}(U')$ and
$\beta$ is in $\psi^{(m+1)N}(U')$.  Any curve complex geodesic
connecting $\beta$ to $\alpha$ will therefore contain a curve $\gamma$
in $\psi^{(m+1)N}(W')$.  Since $\psi^{(m+1)N}(W')$ is contained in
$\psi^{mN}(V')$ the curve $\gamma$ is contained in $\psi^{mN}(V')$, and
we have $d_{\cC}(\beta, \gamma) \ge m + 1$ by induction.  But since $\gamma$ lies on a
geodesic joining $\alpha$ to $\beta$, we have
\begin{eqnarray*}
d_{\cC}(\alpha, \beta) & = & d_{\cC}(\alpha, \gamma) + d_{\cC}(\gamma,
\beta) \\
& \ge & m + 2
\end{eqnarray*}
completing the induction.

Let $V_0 = \cup_{i=0}^{N-1} \psi^i(V')$. There is an $M>0$ such that
$\psi^{MN}(V_0) \subset V'$. Let $n$ be an integer greater than $M$
and assume $k$ is a non-negative integer less than $N$. We then
observe that
$$\psi^{nN+k}(V') \subset \psi^{nN}(V_0) \subset \psi^{n-M}(V')$$
and therefore
$$d_{\cC}(\alpha, \beta) \ge n - M + 1$$
for any $\alpha \in \psi^{nN+k}(V')$ and $\beta \in U'$. Alternatively
if $\alpha \in \psi^n(V')$ and $\beta \in U'$ then
$$d_{\cC}(\alpha, \beta) \ge \left\lfloor \frac{n}{N} \right\rfloor - M + 1 \ge \frac{n}{N} - M.$$

We now set $V'' = \psi^{(B+M)N}(V')$.  Any $\alpha \in V''$ has image
$\psi^n(\alpha)$ lying in $ \psi^{(B+M)N + n}(V')$ so we have
$$d_{\cC}(\psi^n(\alpha), \beta) \ge \frac{(B+M)N + n}{N} - M = \frac{n}{N} + B$$
and $U'$ and $V''$ satisfy (1).

We may assume that $B>D$ and let $W = \pml(S)\backslash (U' \cup V'')$ so
that any path from $U'$ to $V''$ contains a subpath of length at least
$D$. Let $N'= (B+M+1)N$ and let $U = \psi^{-N'}(U')$ and $V =
\psi^{N'}(V'')$. Since $U\subset U'$ and $V \subset V''$ (2) will
still hold $U$, $V$ and $W$. We also note that the sets $U$ and
$\psi^{-2N'}(V)$ will satisfy (1) and $W$ is contained in
$\psi^{-2N'}(V)$. A similar statement holds for $\psi^{2N'}(U)$ and
$V$ with $W$ contained in $\psi^{2N'}(V)$. Therefore (3) will hold.

For (4) we note that we can replace $U$, $V$ and $W$ with $\psi^n(U)$, $\psi^n(V)$ and $\psi^n(W)$
for any integer $n$. We also note that
$$\bigcup_{n\in \integers} \psi^n(V) = \pml(S) \backslash [\mu^+] \ \
\ \mbox{ and } \ \ \  \bigcup_{n \in \integers} \psi^n(U) = \pml(S) \backslash [\mu^-]$$
and that
\begin{eqnarray*}
\bigcup_{n \in \integers} \psi^n(W) & = & \pml(S)\backslash
\left(\bigcap_{n \in \integers} \psi^n(V) \cup \bigcap_{n \in
    \integers} \psi^n(U)\right)  \\ 
&=& \pml(S) \backslash \{[\mu^+], [\mu^-]\}.
\end{eqnarray*}
If we want a fixed curve $\beta$ to be in $U$ we choose $n$ large
enough such that $\beta$ is in $\psi^n(U)$ and then replace $U$ with
$\psi^n(U)$, $V$ with $\psi^n(V)$ and $W$ with $\psi^n(W)$.
\qed{soutoAnosov}

\newcommand{\depth}{\mathbf{depth}}
\newcommand{\virto}{\dashrightarrow}
\newcommand{\evento}{\dashrightarrow}

\subsection{Convergence of iteration}
Let $S$ be a closed surface.  Given $Y \in \Teich(S)$, the {\em Bers
  slice} 
$$B_Y = \{ Q(X,Y) \colon X \in \Teich(S)\} \subset AH(S)$$
parametrizes $\Teich(S)$ by quasi-Fuchsian manifolds.  Since the Bers
slice $B_Y$ has compact closure in $AH(S)$ \cite{Bers:bdry}, the sequence
$\{Q(\psi^{-n} (X), Y)\}_n$ will have a convergent subsequence. It was
established in \cite{McMullen:book:RTM} via a geometric limit argument
that the sequence converges.  We give a new proof as an example of our
methods.

\begin{theorem}
\label{bersconv}
Let $\psi \in \Mod(S)$ be a pseudo-Anosov mapping class. Then the iteration
$\{Q(\psi^n (X), Y)\}_n$ converges in $AH(S)$.
\end{theorem}

\bold{Proof.} Let $T_n$ be the distance between the components of the
boundary of convex cores of $Q(\psi^{n} X, Y)$.  Choose closed
geodesics $\alpha$ on $X$ and $\beta$ on $Y$ of length less than the
Bers constant $L_g$ for $S$.  Then applying
Theorem~\ref{soutoAnosov} there is a positive integer $N$, so that
$$d_\cC(\psi^{n+N}(\alpha),\beta) \ge Kn.$$
Since $\psi^{n+N}(\alpha)$ and $\beta$ have lengths on $\psi^{n+N}(X)$
and $Y$ (respectively) bounded by $L_g$, we have
$$d_{\cC}(\psi^{n+N}(X), Y) = d_{\cC}(\psi^{n+N}(\alpha), \beta) \ge Kn.$$
Applying Corollary \ref{quasifuchsianthick}, we have
$$T_n > K_1 n - K_2.$$
The pseudo-Anosov mapping class $\psi$ is uniformly quasiconformal as
a mapping from $\psi^{n} (X)$ to $\psi^{n+1} (X)$ (independent of
$n$). Let $d_n$ be the distance between $Q(\psi^{n} X, Y)$ and
$Q(\psi^{n+1} X, Y)$ in the Bers slice $B_Y$ as in
Theorem~\ref{qcschwarzinflex}.  Then there are constants $C_1$ and
$C_2$ such that
\begin{eqnarray*}
d_n &\leq& C_1 e^{-C_2 T_n} \\
& \leq & C_1 e^{-C_2(K_1 n - K_2)}.
\end{eqnarray*}
This implies that $Q(\psi^{n}(X), Y)$ is a Cauchy sequence and hence
convergent. \qed{bersconv}

\subsection{Double Limits}
We now examine the pseudo-Anosov double iteration
$$Q_n = Q(\psi^{-n} (X),
\psi^n (Y)).$$ 
Thurston's {\em double limit theorem}
\cite{Thurston:hype2} guarantees that the sequence has a convergent
subsequence.  In \cite{Cannon:Thurston:peano} a proof of convergence
is outlined that uses the Mostow rigidity theorem -- in
\cite{McMullen:book:RTM}, McMullen showed convergence explicitly.

Using the geometric inflexibility theorem proven here, we will give a
single unified proof of this and other convergence convergence results
that is independent of Thurston's double limit theorem.

\begin{theorem}
  \label{doubleconverges} {\sc (Pseudo-Anosov Double Limits)}
  Given $X$, $Y$ in $\Teich(S)$, and a pseudo-Anosov mapping
  class $\psi \in \Mod(S)$, the double iteration $\{Q_n\}_n$ converges
  strongly in $AH(S)$.
\end{theorem}

For each $n$ there is a $K$-quasi-conformal deformation from $Q_n$ to
$Q_{n+1}$ where $K$ bounds the Teichm\"uller distance from $X$ to
$\psi^{-1}(X)$ and from $Y$ to $\psi(Y)$.  Let $\phi_n : Q_n \to Q_{n+1}$
be the map given by Theorem \ref{thickinflex}.  We note that in the
application of the inflexibility theorems, the constants that arise in
Theorem \ref{thickinflex} and \ref{qclengthinflex} depend on $K$ but
not on $n$.

We begin with a criterion to ensure that a curve $\gamma$ has a
convergent sequence of geodeisc lengths $\ell_{Q_n}(\gamma)$.
\begin{prop}
\label{lengthconverges}
Given $L >0$ there exists $B>0$ so that the following holds: if given
$\gamma \in \calC^0(S)$ there exists $N \in \natls$ for which
$\ell_{Q_N}(\gamma) < L$ and for all $n$ we have
$$\min\{d_\calC(\psi^{N+n}(Y), \gamma),
d_\calC(\psi^{-N-n}(X), \gamma)\} \ge K_\psi n +B$$
then there exists $\ell_\infty >0$ so that $\ell_{Q_n}(\gamma) \to
\ell_\infty$.
Furthermore, we have 
$$\ell_{Q_{N+n}}(\gamma) \leq 2
\ell_{Q_N}(\gamma)$$ for all $n>0$.
\end{prop}

\bold{Proof.}  
Let $\depth_Q(\gamma)$ denote the distance of $\tube'_{\varepsilon_3}(\gamma)$
from the 
boundary of the convex core of $Q$.
Let $f_{2L}$ be the function given by Corollary~\ref{curvedeepincore} so that any
curve $\beta$ for which $\ell_{Q}(\beta) < 2L$ satisfies
$$\depth_Q(\beta) \ge f_{2L}( \min\{d_\calC(X,\gamma), d_\calC(\gamma,Y)\}).$$

Let $d_n = \depth_{Q_{N+n}}(\gamma)$.  If $\ell_{Q_{N+n}}(\gamma) <
2L$, then, we have
$$d_{n} \ge  f_{2L}(K n + B).$$
Since $\ell_{Q_N}(\gamma) < L$, we know by assumption that
$$d_{0} \ge f_{2L}(B).$$
Note that $f_{2L}$ is an increasing function so we can make $f_{2L}(B)$ as large as we like through our choice of $B$.   

Let
$$\ell_n = \ell_{Q_{N+n}}(\gamma).$$
If $\ell_n < 2L$ then by Theorem~\ref{qclengthinflex}
\marginpar{right?} there are constants $C'_1$ and $C'_2$ so that
\begin{equation}
\label{logbound}
\left| \log \left(  \frac{\ell_{n+1}}{\ell_n} \right) \right| 
\le  C'_1 e^{ - C'_2 d_n}.
\end{equation}
Choose $C_1$ and $C_2$ such that
$$C_1 e^{-C_2 n} = C'_1 e^{-C_2'f_{2L}(K_\psi n + B)}$$ (recall
$f_{2L}$ is linear).
Since $C_1 = C'_1e^{-C_2' f_{2L}(B)}$ we can choose $B$ such that
$$\sum_{j= 0}^n C_1 e^{-C_2 n} \le
\frac{ C_1 }{1 - e^{-C_2}} \le \log 2$$
from which it follows that
\begin{equation}
\label{log2}
 \left| \log \left(  \frac{\ell_{n}}{\ell_0} \right) \right| 
 < \log 2
\end{equation}
by a simple inductive argument.  We conclude that
$\ell_n < 2L$, and thus equation~\eqref{logbound} holds for all $n$.
It follows that $\{\ell_n\}$ is a Cauchy sequence, and its convergence
to $\ell_\infty >0$
follows from~(\ref{log2}).
\qed{lengthconverges}

We note the following corollary, which will play a role in
establishing strong convergence of $\{Q_n\}$.
\begin{cor}
  If $Q_n$ has a subsequence that converges algebraically then for
  each $\gamma \in \cC^0(S)$ the sequence $\ell_{Q_n}(\gamma)$
  converges to a positive number.
\label{liminf}
\end{cor}

\bold{Proof.}  Algebraic convergence of the subsequence $Q_{n_i}$
implies there is an $L$ so that $\ell_{Q_{n_i}}(\gamma) <L$ for each
$n_i$.  Taking $B$ as guaranteed by Proposition~\ref{lengthconverges}
we use Theorem \ref{soutoAnosov} to choose subsets $U$, $V$ and $W$ of
$\pml(S)$ with $\gamma \in W$ such that
$$\min\{d_\calC(\psi^{N+n}(\beta), \gamma),
d_\calC(\psi^{-N-n}(\alpha), \gamma)\} \ge K_\psi n +B$$
for all $\alpha \in U$ and $\beta \in V$ for all $n \ge 0$. We then choose $N>0$ such that if $\psi^{N}(X) \subset V$ and $\psi^{N}(Y) \subset U$. We then have
$$\min\{d_\calC(\psi^{N+n}(Y), \gamma),
d_\calC(\psi^{-N-n}(X), \gamma)\} \ge K_\psi n +B$$
and by Theorem \ref{lengthconverges} we have that $\ell_{Q_n}(\gamma)$ converges to a positive number. \qed{liminf}

We now use Theorems \ref{soutoAnosov} and \ref{lengthconverges} to
find a pants decomposition whose lengths in $Q_n$ converge.

\begin{prop}\label{pantsconverges}
  There exists a pants decomposition $\cP$ such that for every $\gamma
  \in \cP$ the sequence $\ell_{Q_n}(\gamma)$ converges to a positive
  number.
\end{prop}

\bold{Proof.} By Lemma \ref{coarsepath} and the remark that follows it
there is $D$-coarse path from $\psi^{-n}(X)$ to $\psi^n(Y)$ consisting
of curves in $\cC^0(S)$ whose lengths are at most $L_g$ in
$Q(\psi^{-n}(X), \psi^n(Y))$.  Let $B$ be the constant given by
Proposition \ref{lengthconverges} where $L = L'_g$ is the Bers
constant for a pants decomposition.

As in the proof of Corollary \ref{liminf} we can find subsets $U,V$
and $W$ of $\pml(S)$ chosen with respect to the constants $B+1$ and
$D$ and a positive integer $N$ such that
$$\min\{d_\calC(\psi^{N+n}(Y), \gamma),
d_\calC(\psi^{-N-n}(X), \gamma)\} \ge K_\psi n +B+1$$
for all $n\geq 0$ and $\gamma \in W$.

In our coarse path from $\psi^{-N}(X)$ to $\psi^N(Y)$ consisting of
curves whose length is at most $L_g$ in $Q_N$ there is a curve $\gamma
\in W$.
Let $f:Z \to Q_N$ be a 1-Lipschitz surface realizing $\gamma$. We can
then extend $\gamma$ to a pants decomposition $\cP$ on such that for
all $\gamma' \in \cP$ we have 
$$L'_g > \ell_Z(\gamma') >\ell_{Q_N}(\gamma').$$ 
Since $d_\cC(\gamma, \gamma') \leq 1$ we have
$$\min\{d_\calC(\psi^{N+n}(Y), \gamma'),
d_\calC(\psi^{-N-n}(X), \gamma')\} \ge K_\psi n +B$$
for all $n \geq 0$. Then Theorem \ref{lengthconverges} implies that $\ell_{Q_n}(\gamma')$ converges for all $\gamma' \in \cP$. \qed{pantsconverges}

In the following proposition we will show that we have exponential
decay of the bi-Lipschitz constant on the iterated image of
sufficiently deep and thick subsets. The proof has the same basic
structure as the proof of Proposition \ref{lengthconverges}.
\begin{prop}\label{subsetconverges} Given $\epsilon, R, L,C>0$
  there exist $B, C_1, C_2>0$ such that the following holds. Assume
  that $\cK$ is a subset of $Q_N$ such that $\diam(\cK)<R$,
$\inj_p(\cK) > \epsilon$ for each 
  $p \in \cK$ and $\gamma \in \cC^0(S)$ is represented by a closed
  curve in $\cK$ of length at most $L$ satisfying
$$\min\{d_\calC(\psi^{N+n}(Y), \gamma),
d_\calC(\psi^{-N-n}(X), \gamma\} \ge K_\psi n +B$$ for all $n \geq 0$.
Then we have
$$\log \bilip (\phi_{N+n}, p) \leq C_1 e^{-C_2 n}$$ 
for $p$ in $\phi_{N+n-1} \circ \cdots \circ \phi_{N} \circ f(\cK)$
and
$$\frac{C_1}{1 - e^{-C_2}} < C.$$
\end{prop}

\bold{Proof.} As in the proof of Proposition \ref{lengthconverges}, if
$\calK$ is a subset of $C(Q)$ we let $\depth_Q(\calK)$ be defined by
the distance from $\calK$ to $\bdry C(Q)$.  Let $\Phi_n = \phi_{N+n}
\circ \cdots \phi_N$.

Let $\epsilon_0 = \epsilon e^{-C}$. By \cite{Brooks:Matelski:collars}
there exists an $\epsilon_1$ such that a point within $\epsilon$ of
point of injectivity radius at least $\epsilon_0$ will have
injectivity radius at least $\epsilon_1$.  Let $C'_1$ and $C'_2$ be
the constants given by Theorem \ref{thickinflex} for the thickness
constant $\epsilon_1$. Let $f = f_{Le^C}$ be the linear function given
by Corollary \ref{curvedeepincore}.  
We then define $d(n) = f(K_\psi n + B) - R-\epsilon$. 
Finally we choose $C_1$ and $C_2$ such that
$$C_1e^{-C_2 n} = C'_1 e^{-C'_2 d(n)}.$$
Note that $C_1 = C'_1 e^{-C'_2 d(0)}$ and we can make $d(0)$ as large
as we like through our choice of $B$. Therefore we can choose $B$ such
that
$$\frac{C_1}{1-e^{-C_2}}  < C.$$

With this setup it is now easy to complete the proof of the theorem
via induction. Note that if $p \in \cK$ then we have
$\depth_{Q_N}(B_\epsilon(p)) > d(0)$ where $B_\epsilon(p)$ is the ball
of radius $\epsilon$ centered at $p$. By Theorem \ref{thickinflex} for
all $q \in B_\epsilon(p)$ we have
$$\log \bilip (\phi_N, q) < C_1.$$
In particular for every $p \in \cK$ we
have
$$\log \bilip (\phi_N, p) < C_1$$
and every point in $\phi_N(\cK)$ has injectivity radius at least
$\epsilon e^{-C_1}$.

Let
$$c_n = \sum_{i=0}^n C_1 e^{-C_2 i}$$
and note that
$$c_n < \frac{C_1}{1 - e^{-C_2}} < C.$$

Assume that the theorem holds for all $i$ between $0$ and $n$ and that
the injectivity radius of every point in $\Phi_i(\cK)$ is greater than
$\epsilon e^{-c_i}>\epsilon e^{-C}$.  Note that if $p$ is in $\cK$
then
$$\log \bilip (\Phi_n, p) < c_n < C.$$
It follows that the length of $\Phi_n(\gamma)$ is $< 2L_g e^C$ and for
every point $p$ within $\epsilon$ of $\Phi_n(\cK)$ we have
$\depth_{Q_{N+n}}(p) >d(n)$. We also note that the injectivity radius
at $p$ will be greater than $\epsilon_1$ so we can apply Theorem
\ref{thickinflex} to see that
$$\log \bilip (\phi_{N+n+1}, p) < C_1 e^{-C_2(n+1)}$$
and that at every point in $\phi_{N+n+1} \circ \Phi_n(\cK) =
\Phi_{n+1}(\cK)$ the injectivity radius is at least
$\epsilon^{-c_{n+1}}$. This completes the proof of the induction
hypothesis and the proposition. \qed{subsetconverges}

\begin{prop}\label{surfaceconverges}
  There exists a positive integer $N$, positive constants $C_1$, $C_2 $
  and a 1-Lipschitz homotopy equivalence $f: Z \to Q_N$ so that for
  all points $p$ in $\phi_{N+n-1} \circ \cdots \phi_N \circ f(Z)$ we
  have
$$\log \bilip (\phi_{N+n}, p) < C_1 e^{-C_2 n}.$$
\end{prop}

\bold{Proof.} By Proposition \ref{pantsconverges} there exists a pants
decomposition $\cP$ such that $\ell_{Q_n}(\gamma)$ converges to a
positive number for every $\gamma \in \cP$. In particular there are
constants $L^+>L^->0$ such that $L^+ > \ell_{Q_n}(\gamma) > L^-$ for
all $n$ and $\gamma \in \cP$. For each $n$ let $f_n:Z_n \to Q_n$ be a
1-Lipschtiz homotopy equivalence realizing $\cP$. By the collar lemma,
there exists
$\epsilon>0$ such that any hyperbolic surface with a pants
decomposition whose lengths are between $L^-$ and $L^+$ is
$\epsilon$-thick. In particular the surfaces $Z_n$ are
$\epsilon$-thick. We also note that there is an $R>0$ such that an
$\epsilon$-thick surface has diameter bounded above by $R$.

By Lemma \ref{catthin} there exists an $\epsilon'>0$ such that
$f_n(Z_n)$ is contained in the $\epsilon'$-thick part of $Q_n$. Let
$B$ be the constant given by Proposition \ref{subsetconverges} for the
constants $\epsilon'$, $R$, $L^+$ and $C=2$. (Note that the choice of $2$
is completely arbitrary and could be any number $>1$).  Using Theorem
\ref{soutoAnosov} we can find an integer $N$ such that
$$\min\{d_\calC(\psi^{N+n}(Y), \gamma),
d_\calC(\psi^{-N-n}(X), \gamma\} \ge K_\psi n +B$$ where $\gamma$ is a
curve in $\calP$.\marginpar{\tiny this used to say ``where $\gamma$ is
  curve.'' okay now?}  When then let $f= f_N$ and $Z=Z_N$ and
proposition follows from Proposition
\ref{subsetconverges}. \qed{surfaceconverges}
 
\medskip

We are now ready to prove the convergence of double iteration,
Theorem~\ref{doubleconverges}.
\medskip

\bold{Proof} {\em (Proof of Theorem~\ref{doubleconverges})}.  Let $f:Z
\to Q_N$ be the 1-Lipschtiz surface given by Proposition
\ref{surfaceconverges}. Then the maps $f_n = \phi_{N+n-1} \circ \cdots
\phi_N \circ f$ are $C$-Lipschitz where
$$C = \frac{C_1}{1 - e^{-C_2}}.$$
By Proposition \ref{algebraiclimit} the sequence $\{(f_n,Q_n)\}$ has a
convergent subsequence $\{(f_{n_i},Q_{n_i})\}$ in $AH(S) = AH(Z)$. Let
$\{(f_\infty,Q_\infty)\}$ be the limit. Note that from the proof of
Proposition \ref{algebraiclimit} we can assume that there are lifts
$\tilde{f}_{n_i}$ that converge to $\tilde{f}_\infty$.

Since Corollary \ref{liminf} guarantees the limit has no parabolics,
Theorem \ref{typestrong} implies that the limit is strong. In
particular if we pick a point $p$ in $Z$ and let $p_n = f_n(p)$ then
the sequence $\{(Q_{n_i}, p_{n_i})\}$ will converge geometrically to
$(Q_\infty, p_\infty)$ for some point $p_\infty \in
Q_\infty$. Furthermore if $\cK$ is a compact set with $f_\infty(S)
\subset \cK$ and $g_{n_i}: (\cK, p_\infty) \to (Q_{n_i}, p_{n_i})$ are
approximating maps then $f_{n_i}$ is homotopic to $g_{n_i} \circ
f_\infty$.

We will show that the entire sequence $\{(Q_n, p_n)\}$ converges
geometrically to $(Q_\infty, p_\infty)$. Let $\gamma \in \cC^0(S)$ be
a simple closed curve on $S$ and represent it by a closed curve
$\gamma_\infty$ in $M_\infty$ and let $L=
\ell_{M_\infty}(\gamma_\infty)$. Let $\cK$ be a compact set in
$M_\infty$ and assume that both $p_\infty$ and $\gamma_\infty$ is
contained in $\cK$. To show geometric convergence we need to show that
for any $A>0$ there exists $e^A$-bi-Lipschitz embeddings
$$g_n:(\cK, p_\infty) \to (Q_n, p_n)$$
for $n$ sufficiently large.

Let $B, C_1$ and $C_2$ be the constants given by Proposition
\ref{subsetconverges} with repsect to the constants $\epsilon
e^{-A/2}, Re^{A/2}, Le^{A/2}$ and $A/2$. By Theorem \ref{soutoAnosov}
there exists an $N_A$ such that
$$\min\{d_\calC(\psi^{N_A+n}(Y), \gamma),
d_\calC(\psi^{-N_A-n}(X), \gamma\} \ge K_\psi n +B.$$

Let $\cK'$ be the closed $\epsilon$-neighborhood of $\cK$. By the
strong convergence of the subsequence for $n_i$ sufficiently large
there is a $e^{A/2}$-bi-Lipschitz embedding
$$g_{n_i}: (\cK', p_\infty) \to (M_{n_i}, p_{n_i}).$$
Note that every point in $g_{n_i}(\cK)$ will have injectivity radius
at least $\epsilon e^{-A/2}$, the diameter of $g_{n_i}(\cK)$ will be
at most $Re^{A/2}$ and the length of $g_{n_i}(\gamma_\infty)$ will be
at most $Le^{A/2}$.  Since we can always replace $N_A$ with a larger
integer we can assume $N_A = n_i$ where $n_i$ is part of the
convergent subsequence. We now apply Theorem \ref{subsetconverges} to
$g_{N_A}(\cK)$ which implies that $\phi_{N_A + n - 1} \circ \cdots
\circ \phi_{N_A}$ is $e^{A/2}$-bi-Lipschitz on
$g_{N_A}(\cK)$. Therefore the composition $$g_{N_A + n} = \phi_{N_A + n
  - 1} \circ \cdots \circ \phi_{N_A} \circ g_{N_A}$$ is
$e^A$-bi-Lipschitz. Furthermore $g_n(p_\infty) = p_n$ so we have our
desired bi-Lipschitz embeddings and $\{(Q_n,p_n)\}$ converges
geometrically to $(M_\infty, p_\infty)$.

To see that the sequence also converges algebraically assume that
$\cK$ contains $f_\infty(S)$. We note that $f_{n+1}$ is homotopic to
$\phi_n \circ f_n$ and more generally $f_{n+k}$ is homotopic to
$\phi_{n+k} \circ \cdots \circ \phi_n \circ f_n$. On the subsequence
$\{n_i\}$ we already know that $g_{n_i} \circ f_\infty$ is homotopic
to $f_{n_i}$. By the above fact, the composition $\phi_{N_A+n-1} \circ
\cdots \circ \phi_{N_A} \circ f_{N_A}$ is homotopic to $f_{N_A + n
  -1}$ and in turn homotopic to $g_{N_A + n -1} \circ f_\infty$. If
$g_n$ are approximating maps for $\cK$ whose bi-Lipschitz constant
limits to $1$ then by Lemma \ref{geomtoalg} we have that $( g_n \circ
f_\infty, Q_n) \rightarrow (f_\infty, M_\infty)$. By the above remarks
$( g_n \circ f_\infty, Q_n) \in [(f_n , Q_n)]$ and therefore $[(f_n,
Q_n)]$ converges to $[(f_\infty, M_\infty)]$
algebraically. \qed{doubleconverges}

We conclude with the proof of Theorem~\ref{hyperbolization}.  

\begin{refthm}{Theorem~\ref{hyperbolization}}{\sc (Mapping Torus Hyperbolic)}
  Let $\psi \in \Mod(S)$ be pseudo-Anosov.  Then the mapping torus
  $T_\psi = S \times [0,1]/(x,0) \sim (\psi(x),1)$ admits a complete
  hyperbolic structure.
\end{refthm}

\bold{Proof.}  We note that as $\Mod(S)$ acts diagonally on
quasi-Fuchsian space by re-marking, the manifolds $Q_n$ and
$\psi(Q_n)$ are isometric.  Because we have
$$d(\psi^{-n+1}(X),\psi^{-n}(X)) = d(\psi(X),X) \ \  \ 
\text{and} \ \ \  d(\psi^{n+1}(Y),\psi^n(Y)) = d(\psi(Y), Y),$$ 
and
$$\psi(Q_n) = Q(\psi^{-n +1}(X),\psi^{n+1}(Y)),$$ there is a uniform
$K$ for which $Q_n$ admits a $K$-bi-Lipschitz self-diffeomorphism
$$\Psi_n \colon Q_n \to Q_n$$ so that
$\Psi_n$ is in the homotopy class of $\psi$.

We now use inflexibility and the fact that $(Q_n, p_n)$ converges
geometrically to $(Q_\infty, p_\infty)$ to extract a limiting isometry
$$\Psi_\infty :Q_\infty \longrightarrow Q_\infty$$
in the homotopy class of $\psi$ as a limit directly.  

If $\cK$ is a compact set in $Q_\infty$ containing $p_\infty$ with
geometric limit mappings $g_n : (\cK,p_\infty) \longrightarrow (Q_n,
p_n)$, then $(g_n)^{-1} \circ \Psi_n \circ g_n$ converges up to
subsequence to a uniformly bi-Lipschitz $\Psi_\infty$ where the
conjugating maps are defined.  Observe that since $\depth_{Q_n}(p_n)
\to \infty$ as $n \to \infty$, the compact sets $g_n(\cK)$ are
arbitrarily deep in the convex core of $Q_n$ as $n \to \infty$.
Theorem~\ref{thickinflex} \marginpar{\tiny which theorem precisely?}
then implies that for any $\epsilon>0$ and any $\cK$ the maps
$\Psi_n$ can be taken to be $(1+\epsilon)$-bi-Lipschitz on $g_n(\cK)$
for $n$ sufficiently large.  Diagonalizing, the limit $\Psi_\infty$ is
an isometry.  As the group of isometries of any hyperbolic 3-manifold
is discrete, we may pass to the quotient by the action of $\langle
\Psi_\infty\rangle$.  The quotient of $Q_\infty/\langle
\Psi_\infty\rangle$ is a hyperbolic 3-manifold with the fundamental
group $\pi_1(T_\psi)$, and is thus homeomorphic $T_\psi$ by Stallings'
Theorem \cite{Stallings:fibers}.
% At this stage, one may simply follow Thurston's
% original approach and use the compactness of quasiconformal maps. In
% particular, the maps $\Psi_n$ lift to $\hthree$ and extend to
% $K'$-quasiconformal maps
% $$\tilde{\psi}_n: \chat \longrightarrow \chat$$
% where $K'$ only depends on $K$. After passing to a subsequence these
% maps will limit to a quasiconformal map $\tilde{\psi}_\infty$. Note
% that if $Q_n = \hthree/\Gamma_n$ then $\Gamma_n$ and $\tilde{\psi}_n$
% will generate a group of homeomorphisms of $\chat$ isomorphic to
% $\pi_1(T_\psi)$ and therefore this property will hold in the
% limit. The Beltrami differential of $\tilde{\psi}_\infty$ will be
% $\Gamma_\infty$-invariant and since the limit set of $\Gamma_\infty$
% is all of $\chat$, Sullivan's Rigidity
% Theorem \cite{Sullivan:linefield} implies that 
% $\tilde{\psi}_\infty$ is conformal. In particular the group
% $\hat{\Gamma}$ generated by $\Gamma_\infty$ and $\tilde{\psi}_\infty$
% acts conformally on $\chat$. It is not hard to show the the
% discreteness of $\Gamma_\infty$ in $\PSL_2(\cx)$ implies that
% $\hat{\Gamma}$ is discrete and therefore we have a discrete subgroup
% of $\PSL_2(\cx)$ that is isomorphic to $\pi_1(T_\psi)$. By a theorem
% of Stallings \cite{Stallings:fibers} we have that
% $\hthree/\hat{\Gamma}$ is homeomorphic to $T_\psi$.
\qed{hyperbolization}

\bold{Remark.}  Note that in Thurston's original proof 
significant extra work is required to show that the limit $Q_\infty$
is doubly degenerate, or that the limit set of $\Gamma_\infty $ is the
entire sphere (see \cite[Sec.~6.2]{Otal:book:fibered}).  In our
approach the double degeneracy of $Q_\infty$ is immediate from our
estimates on the depth of the basepoint in the convex core and the
strong convergence of $Q_n$.

\bibliographystyle{math}
\bibliography{math}

\begin{thebibliography}{CCHS}

\bibitem[Ag]{Agol:tame}
I.~Agol.
\newblock {Tameness of hyperbolic 3-manifolds}.
\newblock {\em Preprint, {\tt arXiv:mathGT/0405568}} (2004).

\bibitem[Ah]{Ahlfors:visual}
L.~Ahlfors.
\newblock {Invariant operators and integral representations in hyperbolic
  space}.
\newblock {\em Math. Scan.} {\bf 36}(1975), 27--43.

\bibitem[And]{Anderson:dirichlet}
M.~Anderson.
\newblock {The Dirichlet problem at infinity for manifolds of negative
  curvature}.
\newblock {\em J. Diff. Geom.} {\bf 18}(1983), 701--721.

\bibitem[Brs1]{Bers:simunif}
L.~Bers.
\newblock {Simultaneous uniformization}.
\newblock {\em Bull. AMS} {\bf 66}(1960), 94--97.

\bibitem[Brs2]{Bers:bdry}
L.~Bers.
\newblock {On boundaries of Teichm\"uller spaces and on kleinian groups: I}.
\newblock {\em Annals of Math.} {\bf 91}(1970), 570--600.

\bibitem[BeFu]{Bestvina:Fujiwara:bounded}
M.~Bestvina and K.~Fujiwara.
\newblock {Bounded cohomology of subgroups of mapping class groups}.
\newblock {\em Geometry and Topology} {\bf 6}(2002), 69--89.

\bibitem[Bon1]{Bonahon:tame}
F.~Bonahon.
\newblock {Bouts des vari\'et\'es hyperboliques de dimension 3}.
\newblock {\em Annals of Math.} {\bf 124}(1986), 71--158.

\bibitem[Bon2]{Bonahon:currents}
F.~Bonahon.
\newblock {The geometry of Teichm\"uller space via geodesic currents}.
\newblock {\em Invent. math.} {\bf 92}(1988), 139--162.

\bibitem[Bow]{Bowditch:hulls}
B.~Bowditch.
\newblock {Some results on the geometry of convex hulls in manifolds of pinched
  negative curvature}.
\newblock {\em Comment. Math. Helv.} {\bf 69}(1994), 49--81.

\bibitem[BH]{Bridson:Haefliger:npc}
M.~Bridson and A.~Haefliger.
\newblock {\em Metric Spaces of Non-Positive Curvature}.
\newblock Springer-Verlag, 1999.

\bibitem[BB]{Brock:Bromberg:density}
J.~Brock and K.~Bromberg.
\newblock {On the density of geometrically finite Kleinian groups}.
\newblock {\em Acta Math.} {\bf 192}(2004), 33--93.

\bibitem[BBES]{BBES:elc}
J.~Brock, K.~Bromberg, R.~Evans, and J.~Souto.
\newblock {Maximal cusps, ending laminations and the classification of Kleinian
  groups}.
\newblock {\em In preparation (2008)}.

\bibitem[BCM]{Brock:Canary:Minsky:elc}
J.~Brock, R.~Canary, and Y.~Minsky.
\newblock {The classification of Kleinian surface groups, II: the ending
  lamination conjecture}.
\newblock {\em Submitted (2004).}

\bibitem[Brm1]{Bromberg:Schwarzian}
K.~Bromberg.
\newblock {Hyperbolic cone manifolds, short geodesics, and Schwarzian
  derivatives}.
\newblock {\em J. Amer. Math. Soc.} {\bf 17}(2004), 783--826.

\bibitem[Brm2]{Bromberg:bers}
K.~Bromberg.
\newblock {Projective structures with degenerate holonomy and the Bers density
  conjecture}.
\newblock {\em Annals of Math.} {\bf 166}(2007), 77--93.

\bibitem[BS]{Bromberg:Souto:density}
K.~Bromberg and J.~Souto.
\newblock {Density of Kleinian groups}.
\newblock {\em In preparation}.

\bibitem[BM]{Brooks:Matelski:collars}
R.~Brooks and J.~P. Matelski.
\newblock {Collars for Kleinian Groups}.
\newblock {\em Duke Math. J.} {\bf 49}(1982), 163--182.

\bibitem[Bus]{Buser:book:spectra}
P.~Buser.
\newblock {\em Geometry and Spectra of Compact Riemann Surfaces}.
\newblock Birkhauser Boston, 1992.

\bibitem[CG]{Calegari:Gabai:tame}
D.~Calegari and D.~Gabai.
\newblock {Shrinkwrapping and the taming of hyperbolic 3-manifolds}.
\newblock {\em J. AMS} {\bf 19}(2006), 385--446.

\bibitem[Can]{Canary:ends}
R.~D. Canary.
\newblock {Ends of hyperbolic 3-manifolds}.
\newblock {\em J. Amer. Math. Soc.} {\bf 6}(1993), 1--35.

\bibitem[CCHS]{Canary:Culler:Hersonsky:Shalen}
R.~D. Canary, M.~Culler, S.~Hersonsky, and P.~B. Shalen.
\newblock {Approximation by maximal cusps in the boundaries of quasiconformal
  deformation spaces}.
\newblock {\em J. Diff. Geom.} {\bf 64}(2003), 57--109.

\bibitem[CH]{Canary:Hersonsky:cusps}
R.~D. Canary and S.~Hersonsky.
\newblock {Ubiquity of geometric finiteness in boundaries of deformation spaces
  of hyperbolic 3-manifolds}.
\newblock {\em Amer. J. of Math.} {\bf 126}(2004), 1193--1220.

\bibitem[CT]{Cannon:Thurston:peano}
J.~W. Cannon and W.~P. Thurston.
\newblock {Group invariant Peano curves}.
\newblock {\em Geometry and Topology} {\bf 11}(2007), 1315--1355.

\bibitem[Dum]{Dumas:handbook}
D.~Dumas.
\newblock {Complex projective structures}.
\newblock In {\em Handbook of Teichm\"uller Theory, Volume II}. EMS Publishing
  House, 2008.

\bibitem[Ev]{Evans:persists}
R.~Evans.
\newblock {Tameness persists in weakly type-preserving strong limits}.
\newblock {\em Amer. J. Math.} {\bf 126}(2004), 713--737.

\bibitem[FLP]{FLP}
A.~Fathi, F.~Laudenbach, and V.~Po\'enaru.
\newblock {\em Travaux de Thurston sur les surfaces}, volume 66-67.
\newblock Ast\'erisque, 1979.

\bibitem[HK]{Hodgson:Kerckhoff:rigidity}
C.~Hodgson and S.~Kerckhoff.
\newblock {Rigidity of hyperbolic cone-manifolds and hyperbolic Dehn surgery}.
\newblock {\em J. Diff. Geom.} {\bf 48}(1998), 1--59.

\bibitem[IT]{Imayoshi:Taniguchi:book}
Y.~Imayoshi and M.~Taniguchi.
\newblock {\em An Introduction to Teichm\"uller Spaces}.
\newblock Springer-Verlag, 1992.

\bibitem[MM]{Masur:Minsky:CCI}
H.~Masur and Y.~Minsky.
\newblock {Geometry of the complex of curves I: hyperbolicity}.
\newblock {\em Invent. Math.} {\bf 138}(1999), 103--149.

\bibitem[Mc1]{McMullen:iter}
C.~McMullen.
\newblock {Iteration on Teichm\"uller space}.
\newblock {\em Invent. math.} {\bf 99}(1990), 425--454.

\bibitem[Mc2]{McMullen:cusps}
C.~McMullen.
\newblock {Cusps are dense}.
\newblock {\em Annals of Math.} {\bf 133}(1991), 217--247.

\bibitem[Mc3]{McMullen:book:RTM}
C.~McMullen.
\newblock {\em Renormalization and 3-Manifolds Which Fiber Over the Circle}.
\newblock Annals of Math. Studies 142, Princeton University Press, 1996.

\bibitem[Min1]{Minsky:torus}
Y.~Minsky.
\newblock {The classification of punctured torus groups}.
\newblock {\em Annals of Math.} {\bf 149}(1999), 559--626.

\bibitem[Min2]{Minsky:CKGI}
Y.~Minsky.
\newblock {The classification of Kleinian surface groups I: models and bounds}.
\newblock {\em Annals of Math.} {\bf 171}(2010), 1--107.

\bibitem[Mum]{Mumford:compact}
D.~Mumford.
\newblock {A remark on Mahler's compactness theorem}.
\newblock {\em Proc. AMS} {\bf 28}(1971), 289--294.

\bibitem[Neh]{Nehari:schwarzian}
Z.~Nehari.
\newblock {Schwarzian derivatives and schlicht functions}.
\newblock {\em Bull. AMS} {\bf 55}(1949), 545--551.

\bibitem[Ot]{Otal:book:fibered}
J.~P. Otal.
\newblock {\em Le th\'eor\`eme d'hyperbolisation pour les vari\'et\'es
  fibr\'ees de dimension trois}.
\newblock Ast\'erisque, 1996.

\bibitem[Rei]{Reimann:visual}
H.M. Reimann.
\newblock {Invariant extension of quasiconformal deformations}.
\newblock {\em Ann. Acad. Sci. Fen.} {\bf 10}(1985), 477--492.

\bibitem[St]{Stallings:fibers}
J.~Stallings.
\newblock {On fibering certain 3-manifolds}.
\newblock In {\em Topology of 3-manifolds}, pages 95--100. Prentice Hall, 1962.

\bibitem[Sul]{Sullivan:linefield}
D.~Sullivan.
\newblock {On the ergodic theory at infinity of an arbitrary discrete group of
  hyperbolic motions}.
\newblock In {\em Riemann Surfaces and Related Topics: Proceedings of the 1978
  Stony Brook Conference}. Annals of Math. Studies 97, Princeton, 1981.

\bibitem[Th1]{Thurston:book:GTTM}
W.~P. Thurston.
\newblock {\em Geometry and Topology of Three-Manifolds}.
\newblock Princeton lecture notes, 1979.

\bibitem[Th2]{Thurston:hype2}
W.~P. Thurston.
\newblock {Hyperbolic structures on 3-manifolds II: Surface groups and
  3-manifolds which fiber over the circle}.
\newblock {\em Preprint, {\tt arXiv:math.GT/9801045}} (1986).

\end{thebibliography}

{\sc \small
 \bigskip

\noindent Brown University \bigskip

\noindent University of  Utah

}

\end{document}